\documentclass[]{scrartcl}

\usepackage[utf8]{inputenc}
\usepackage{amssymb}
\usepackage{graphicx}
\usepackage{multirow}
\usepackage{lscape}
\usepackage{longtable}
\usepackage{xifthen}
\usepackage{lineno}
\usepackage{float}
\usepackage{graphicx}
\graphicspath{ {./src/img/} }

\usepackage{natbib}

\usepackage[dvipsnames]{xcolor}

\usepackage{amsmath}

\usepackage{tikz}
\usetikzlibrary{positioning,graphs,matrix}

\usepackage[hidelinks]{hyperref}

\title{Marine Autonomous Vehicle Fleet Scheduling to Maximise Scientific Impact}
\author{Mehdi El Krari\thanks{Corresponding author: elkrari@acm.org}, Jonathan Smith, Maria Fox \\ British Antarctic Survey, AI Lab\\ Cambridge, UK}

\allowdisplaybreaks
\begin{document}
\newcommand{\duration}[1]{\textrm{\textbf{d}}_{#1}}
\newcommand{\durationGT}[2]{\textrm{\textbf{d}}_{#1,#2}} 
\newcommand{\durationCF}[2]{\textrm{\textbf{d}}_{#1,#2}}
\newcommand{\length}[1]{\textrm{\textbf{l}}_{#1}}
\newcommand{\pTime}{\textit{\textbf{t}}}
\newcommand{\maximum}[1]{\textrm{\textbf{max}}_{#1}}
\newcommand{\minimum}[1]{\textrm{\textbf{min}}_{#1}}
\newcommand{\pNumber}{\textrm{nb}}
\newcommand{\card}[1]{|#1|}
\newcommand{\cost}[1]{\textrm{\textbf{ec}}_{#1}}
\newcommand{\gain}[1]{\textrm{\textbf{sg}}_{#1}}
\newcommand{\capacity}[3]{\textrm{\textbf{cp}}_{#1,#2,#3}}
\newcommand{\requirement}[1]{\textrm{\textbf{rq}}_{#1}}
\newcommand{\requirementType}[2]{\textrm{\textbf{rq}}_{#1,#2}} 
\newcommand{\surveyDuration}[1]{\textrm{\textbf{sd}}_{#1}}

\newcommand{\pSetAUVType}{\textit{VT}}
\newcommand{\setAUVType}[1][]{\pSetAUVType\ifthenelse{\isempty{#1}}{}{_{#1}}}
\newcommand{\AUVType}{\MakeLowercase{\pSetAUVType}}
\newcommand{\numberAUVTypes}{\card{\setAUVType}}
\newcommand{\getAUVType}[1]{\textrm{get}\setAUVType_{#1}}
\newcommand{\pSetAUV}{\textit{MV}}
\newcommand{\setAUV}[1][]{\pSetAUV\ifthenelse{\isempty{#1}}{}{_#1}}
\newcommand{\AUV}[1][]{\MakeLowercase{\pSetAUV}\ifthenelse{\isempty{#1}}{}{_#1}}
\newcommand{\numberAUVs}{\card{\setAUV}}
\newcommand{\getSetAUV}[1]{\textrm{get}\setAUV_{#1}}
\newcommand{\pSetRoute}{\textit{RT}}
\newcommand{\setRoute}[1][]{\pSetRoute\ifthenelse{\isempty{#1}}{}{_#1}}
\newcommand{\route}[2]{\MakeLowercase{\pSetRoute}\ifthenelse{\isempty{#1}}{}{_{#1,#2}}}
\newcommand{\numberRoutes}{\card{\pSetRoute}}
\newcommand{\subSetRoute}[2]{\textit{RT}_{#1,#2}}
\newcommand{\pSetStation}{\textit{RS}}
\newcommand{\setStation}[1][]{\pSetStation\ifthenelse{\isempty{#1}}{}{_{#1}}}
\newcommand{\station}{\MakeLowercase{\pSetStation}}
\newcommand{\numberStations}{\card{\setStation}}
\newcommand{\setBase}{\textit{BS}}
\newcommand{\base}[1]{\MakeLowercase{\setBase}_{#1}}
\newcommand{\setPickup}{\textit{PL}}
\newcommand{\pickup}[1]{\MakeLowercase{\setPickup}_{#1}}
\newcommand{\setLocation}{\textit{LC}}
\newcommand{\location}{\MakeLowercase{\setLocation}}
\newcommand{\setPeriod}{\textit{PR}}
\newcommand{\period}{\MakeLowercase{\setPeriod}}
\newcommand{\pSetTaskReq}{\textit{TR}}
\newcommand{\setTaskReq}[2]{\pSetTaskReq\ifthenelse{\isempty{#1}}{}{_{#1,#2}}}
\newcommand{\taskReq}[2]{\MakeLowercase{\pSetTaskReq}\ifthenelse{\isempty{#1}}{}{_{#1,#2}}}
\newcommand{\taskReqWeight}[1]{\textbf{w}_{#1}}

\newcommand{\batteryLevel}[1]{\textrm{\textbf{bl}}_{#1}}
\newcommand{\initBatteryLevel}[1]{\batteryLevel{#1}^{i}}
\newcommand{\activeBatteryCons}[1]{\textrm{\textbf{bc}}^{a}_{#1}}
\newcommand{\idleBatteryCons}[1]{\textrm{\textbf{bc}}^{i}_{#1}}
\newcommand{\batteryLevelThreshold}[1]{\textrm{\textbf{th}}_{#1}}
\newcommand{\energyFactor}[1][]{\textrm{\textbf{ef}}\ifthenelse{\isempty{#1}}{}{_{#1}}}
\newcommand{\energyFactorCF}[2]{\textrm{\textbf{ef}}_{#1,#2}}
\newcommand{\setShip}{\textit{SH}}
\newcommand{\ship}[1][]{\MakeLowercase{\setShip}\ifthenelse{\isempty{#1}}{}{_{#1}}}

\newcommand{\setTimeWindow}{\textit{TW}}
\newcommand{\timeWindow}[2]{\MakeLowercase{\setTimeWindow}\ifthenelse{\isempty{#1}}{}{_{#1,#2}}}
\newcommand{\earliestDeployment}[1]{\textbf{ed}_{#1}}
\newcommand{\latestDeployment}[1]{\textbf{ld}_{#1}}
\newcommand{\earliestRecovery}[1]{\textbf{er}_{#1}}
\newcommand{\latestRecovery}[1]{\textbf{lr}_{#1}}
\newcommand{\startTime}[1]{\textrm{\textbf{st}}_{#1}}
\newcommand{\eendTime}[1]{\textrm{\textbf{et}}_{#1}}
\newcommand{\earliestTime}[1]{\textrm{\textbf{rt}}_{#1}}
\newcommand{\latestTime}[1]{\textrm{\textbf{lt}}_{#1}}

\newcommand{\dataAUV}[1]{\data_{\AUV_#1}}

\newcommand{\numberDays}{\textit{N}}
\newcommand{\numberPeriods}{\textit{M}}
\newcommand{\deadline}{\textit{DL}}
\newcommand{\largeConstant}{\textrm{K}}

\newcommand{\dAUVStation}[2]{\textrm{x}_{#1,#2}}
\newcommand{\dAUVRoute}[2]{\textrm{y}_{#1,#2}}
\newcommand{\dStationST}[1]{\textrm{bs}_{#1}}
\newcommand{\dAUVStationST}[2]{\textrm{stt}_{#1,#2}}
\newcommand{\dAUVRouteAT}[2]{\textrm{atr}_{#1,#2}}
\newcommand{\dArrivalTime}[3]{\textrm{at}_{#1,#2,#3}}

\newcommand{\dBatteryLevelTask}[2]{\textrm{blt}_{#1,#2}}
\newcommand{\dBatteryLevelRoute}[2]{\textrm{blr}_{#1,#2}}
\newcommand{\dAUVChargingFromTo}[3]{\textrm{z}_{#1,#2,#3}}
\newcommand{\dAUVChargingIn}[2]{\textrm{zi}_{#1,#2}}
\newcommand{\dAUVChargingOut}[2]{\textrm{zo}_{#1,#2}}
\newcommand{\dAUVChargingInRoute}[3]{\textrm{ziy}_{#1,#2,#3}}
\newcommand{\dAUVChargingOutRoute}[3]{\textrm{zoy}_{#1,#2,#3}}

\newcommand{\enumer}[4][,]{\{ #4_{1}#1 #4_{2} #1 \dots #1 #4_{#3-1} #1 #4_{#3} \} = #2}
\flushbottom
\maketitle
\begin{abstract}
\textbf{Abstract:}\\
The marine science community increasingly relies on Marine Autonomous Vehicles (MAVs) to collect the critical environmental data required to understand global ocean systems. However, as these operations scale, manually routing and planning large autonomous fleets becomes exponentially complex and time-consuming. To address this, we propose a mixed-integer linear programming (MILP) model designed to automate and optimise MAV deployment schedules. The model accounts for strict operational constraints, including battery capacities and time windows for data collection, while aiming to maximise total data collection and minimise both the number of deployed vehicles and their energy consumption. A key novelty of this framework is integrating conventional ship itineraries, allowing MAVs to support vessels with mid-mission battery swapping or accelerated transit between waypoints. Computational experiments demonstrate that the model is highly scalable, solving routing problems for fleets of dozens of MAVs in seconds, and scaling to hundreds of vehicles in only a few minutes. Beyond operational scheduling, the framework serves as a robust simulation tool for evaluating 'what-if' scenarios and analysing the impact of varying parameters on deployment strategies. Finally, the solution generates a suite of visualisations designed to enhance explainability and support strategic decision-making for stakeholders. 
\end{abstract}

\clearpage

\section{Introduction}

Marine science is inextricably linked to global environmental health, providing the foundation for understanding and mitigating anthropogenic impacts on Earth's ocean systems. Developing effective protection strategies, deciphering the profound effects of climate change, and driving higher-resolution environmental models all rely on the continuous acquisition of high-fidelity marine data (e.g temperature, salinity, and conductivity). This need is particularly urgent in highly sensitive environments like the polar regions, where rapid climatic shifts have global consequences. To address this escalating demand for continuous ocean observations, the British government—through UK Research and Innovation (UKRI) and its Natural Environment Research Council (NERC)—established the Future Marine Research Infrastructure (FMRI) programme. This strategic initiative is designed to deliver \textit{``the next generation of large-scale marine research infrastructure}\footnote{The Future Marine Research Infrastructure programme: https://fmri.ac.uk/}''.

Marine Autonomous Vehicles (MAVs) \citep{Alam2014, Schofield2007} are integral to this future infrastructure. Driven by the UKRI Environmental Sustainability Strategy to achieve net-zero by 2040, these vehicles offer vastly improved carbon efficiency compared to conventional research vessels, making them increasingly prominent in data collection efforts. Consequently, the FMRI programme is driving a significant expansion of the national autonomous fleet. However, as the number of deployed MAVs grows, the manual planning of their missions becomes exponentially more complex. This process cannot scale effectively without the integration of advanced, automated decision-support tools.

To address this operational challenge, this paper proposes a mixed integer linear programming (MILP) model to optimise routing and deployment planning for large MAV fleets. The developed model accounts for numerous asset- and environment-specific operational constraints to ensure that every generated solution is practically feasible, with the primary objective of maximising the total volume of scientific data collected. Furthermore, our proposed framework is highly adaptable. It integrates pre-existing ship itineraries, allowing conventional vessels to be used for MAV battery swapping or as transit corridors.

Our work can also serve as a simulation tool: researchers can evaluate "what-if" scenarios (e.g., varying battery capacities or extended deployment periods) prior to physical deployment. Finally, to improve the explainability of these complex plans, we introduce a suite of textual and visual outputs derived from the model's solutions, specifically designed to highlight the most critical aspects of the overarching deployment strategy.

The remainder of this paper is structured as follows. Section \ref{sec:literature} presents a review of related literature. Section \ref{sec:formulation} formally defines the problem, followed by the presentation of the mathematical model in Section \ref{sec:model}. Section \ref{sec:case_studies} introduces three distinct use cases alongside their respective computational results. Finally, Section \ref{sec:conclusion} concludes the paper by discussing our findings and perspectives on future work.


\section{Literature Review}
\label{sec:literature}

Marine Autonomous Vehicles (MAVs) are indispensable tools for marine science and naval operations, enabling data collection in environments where human presence is impractical or hazardous \citep{Bovio2006}. Their navigation and operational systems face unique challenges, including limited underwater communication, dynamic ocean currents, and strict energy constraints \citep{Zhang2023}. While recent reviews highlight advancements in sensor fusion and path planning (e.g., \citet{Zhang2023}), operational reliability and fleet-level coordination remain critical bottlenecks for large-scale mission success. Human-in-the-loop frameworks, such as those by \citet{Somers2016}, have traditionally enhanced navigational robustness by combining operator expertise with probabilistic planning, but these manual methods struggle to scale as fleet sizes increase.

\subsection{MAV Path and Mission Planning}

Path planning is fundamental to MAV operations, and early work by \citet{Warren1990} established route-planning techniques for underwater environments. Modern approaches address increasingly complex scenarios, as \citet{Li2018} comprehensively reviews and categorises methods into sampling-based, optimisation-based, bio-inspired, and geometric model search algorithms.

Optimisation-based approaches are particularly relevant for fleet coordination. \citet{Petres2007} demonstrated efficient path planning using Fast Marching methods \citep{sethian1999fast}, while \citet{Yang2022} proposed time-saving solutions for complex underwater conditions, and \citet{Yan2022} introduced a whale-based optimisation algorithm. A critical challenge in underwater path planning is handling environmental uncertainty, which \citet{Yao2021} addressed through interval optimisation for uncertain flow fields.

Mission planning for MAVs involves both high-level task allocation and low-level path planning, which is particularly challenging in dynamic ocean environments. \citet{Tavana2010} proposed a decision-support system that combines multi-criteria analysis and fuzzy sets to address uncertainty in MAV operations. In contrast, \citet{Yan2014} developed an integrated algorithm for mission assignment and path planning in ocean current environments. Evolutionary algorithms have also shown promise for mission planning, as evidenced by \citet{Geng2013}'s work on UAV surveillance, which decomposes the problem into vantage point selection and path optimisation.

\subsection{MAV Fleet Coordination}

While individual MAVs are pivotal for marine data collection, scaling deployments to large fleets introduces immense challenges in mission planning \citep{Sahoo2019, Thompson2019}. Early operational research focused on single-MAV navigation \citep{Leonard2016, Miller2010}, while more recent studies have begun to address multi-MAV coordination. For instance, \citet{Sorbi2012} proposed behaviour-based planners for cooperative MAVs, and \citet{McMahon2016} optimised missions in dynamic environments.

However, historically, the literature has treated MAV deployments as isolated, small-scale scientific endeavours rather than integrated components of a broader maritime logistics network. Consequently, a significant gap remains: existing autonomous coordination models almost universally assume MAVs operate independently of conventional maritime infrastructure. They rarely account for ship support, such as mid-mission battery swaps or assisted transits. This omission stems from the traditional divide between micro-level robotic path planning and macro-scale maritime logistics. By bridging these two fields, our proposed MILP model directly addresses this critical operational reality and reflects the true logistical demands of modern marine research.

\subsection{Optimisation Models in Marine Environments}

Maritime logistics has seen significant advances in optimisation techniques, particularly for fleet management and routing problems. \citet{Brouer2014} established a foundational MILP model for liner shipping network design, demonstrating how integer programming can optimise vessel routes while handling complex operational constraints. Their work with real-world data proved the practical applicability of such models at an industrial scale. For fuel supply vessels, \citet{Christiansen2017} compared arc-flow and path-flow formulations, showing how to model time-dependent constraints effectively. At the strategic level, \citet{Li2021} developed a comprehensive MILP model for global cargo flows that incorporates both containerised and bulk shipping.

Bridging the gap to autonomous fleets, MILP modelling has been applied to MAV path planning under energy constraints \citep{DeCarolis2017}. \citet{Yilmaz2009} worked on a MILP model for adaptive sampling, while \citet{Eichhorn2015} optimised routes in time-varying currents. These studies, however, neglect fleet-scale optimisation and joint ship-MAV planning. Recent work continues to advance MILP applications for autonomous fleets; \citet{Wang2024} demonstrated how constraint programming can solve task allocation as a Capacitated Vehicle Routing Problem (CVRP), and \citet{Pang2021} showed quantum particle swarm optimisation's effectiveness for energy-constrained path planning.

\subsection{Energy Constraints and Fleet Scalability}

Battery management is the primary limiting factor for persistent MAV operations. \citet{DeCarolis2017} estimated runtime energy, while \citet{Bychkov2019} and \citet{MahmoudZadeh2018Bis} optimised rendezvous strategies with stationary recharge stations. \citet{Yazdani2016} further demonstrated the importance of battery and energy management through docking operations, developing a real-time trajectory optimiser that minimises both mission time and energy expenditure while accounting for hydrodynamic effects.

While evolutionary algorithms \citep{MahmoudZadeh2018} and particle swarm methods \citep{Lu2020} can handle energy-aware task allocation, they lack the rigorous mathematical guarantees of MILP for optimality. Furthermore, prior research has largely avoided exact MILP formulations for joint ship-MAV routing due to the exponential computational complexity introduced by combining continuous autonomous operations with discrete maritime vessel schedules.

Operating within the context of next-generation, large-scale marine research infrastructure, our unique position allows us to formulate a highly scalable solution that mirrors actual deployment constraints. Our model bridges the gap between theoretical multi-agent routing and applied maritime logistics, ensuring strict operational feasibility for hundreds of MAVs while explicitly incorporating moving support vessels—a significant scalability and logistical leap over prior heuristic approaches \citep{Wang2022}.
\section{Problem Formulation}
\label{sec:formulation}

This work distinguishes itself by tackling the joint routing of autonomous fleets and conventional support vessels—a problem previously avoided because of its inherent complexity. Rather than relying on stationary infrastructure or heuristic approximations, our framework introduces a highly scalable MILP model that explicitly incorporates discrete, scheduled ship itineraries. This approach enables the dynamic assignment of MAVs to support vessels for mid-mission battery swapping and accelerated transit, directly addressing the operational realities of large-scale marine research.

\begin{sloppypar}
We define a problem consisting of $\numberStations$ research stations, denoted by the set $\enumer{\setStation}{\numberStations}{\station}$, and $\numberAUVs$ Marine Autonomous Vehicles (MAVs) denoted by the set $\enumer{\setAUV}{\numberAUVs}{\AUV}$. Each MAV $\AUV \in \setAUV$ possesses a unique type $\AUVType \in \enumer{\setAUVType{}}{\numberAUVTypes}{\AUVType}$ that defines its ability to perform some task requirements. 

The set of MAVs of a given type $\AUVType$ is $\setAUV[\AUVType]$, $\enumer[\cup]{\setAUV}{\numberAUVTypes}{{\setAUV[\AUVType]}}$ and $\enumer[\cap]{\varnothing}{\numberAUVTypes}{{\setAUV[\AUVType]}}$.
\end{sloppypar}

Each MAV $\AUV \in \setAUV$, if assigned to at least one task requirement, should be deployed from the base, known as port of mobilisation, $\base{\AUV}$ within a time window $\timeWindow{\earliestDeployment{\base{\AUV}}}{\latestDeployment{\base{\AUV}}}$, and be recovered at its pickup location (or port of demobilisation) within the time window $\timeWindow{\earliestRecovery{\pickup{\AUV}}}{\latestRecovery{\pickup{\AUV}}}$. $\earliestDeployment{\base{\AUV}}, \latestDeployment{\base{\AUV}}, \earliestRecovery{\pickup{\AUV}}, \latestRecovery{\pickup{\AUV}}$ are respectively the earliest and latest time for deployment from the base $\base{\AUV}$ of $\AUV$, and the earliest and latest time for recovery at the pickup location $\pickup{\AUV}$ of $\AUV$. 

\subsection{Task Requirements}
The problem has a set of task requirements $\setTaskReq{}{}$; every single task $\taskReq{\station}{\AUVType} \in \setTaskReq{}{}$ maps a research station $\station$ with an MAV type $\AUVType$ that can cover the task. These requirements are defined by the following characteristics:
\begin{itemize}
    \item $\requirement{\taskReq{}{}}$: the required number of MAVs of type $\AUVType$ needed for completing the task requirement $\taskReq{}{}$. The number of assigned MAVs to each task $\taskReq{}{}$ can be less or equal than $\requirement{\taskReq{}{}}$ in the optimal solution, but not greater;
    \item $\duration{\taskReq{}{}}$: the duration for the task $\taskReq{}{}$ to be completed. Each assigned $\AUV$ to $\taskReq{}{}$ should spend exactly $\duration{\taskReq{}{}}$ without interruption;
    \item $\earliestTime{\taskReq{}{}}$: the earliest start time for the task $\taskReq{}{}$; 
    \item $\latestTime{\taskReq{}{}}$: the latest time for finishing the task $\taskReq{}{}$;
    \item $\taskReqWeight{\taskReq{}{}}$: the scientific weight associated to the task requirement $\taskReq{}{}$. $\taskReqWeight{\taskReq{}{}}$ is included in the objective function (\ref{eq:set_objective}) to prioritise higher-weighted scientific missions over others under resource and vessel availability constraints.
\end{itemize}

\subsection{Battery constraints} 

Each MAV $\AUV \in \setAUV$ has a battery level $\batteryLevel{\AUV}$ that varies every unit of time, from deployment, starting with its initial battery level $\initBatteryLevel{\AUV}$ and getting lower until its recovery. Battery usage depends on the MAV type and environmental conditions. 

Each MAV type $\AUVType \in \setAUVType{}$ has two different average battery consumptions per time unit, either when it is active $\activeBatteryCons{\AUVType}$ (i.e. doing review tasks) or idle $\idleBatteryCons{\AUVType}$ (e.g. transiting). 
In addition to the consumptions, $\AUVType$ also defines the minimum battery level (threshold) $\batteryLevelThreshold{\AUVType}$ for each vehicle $\AUV \in \setAUV[\AUVType]$.

Depending on its location, each MAV's battery consumption can be higher or lower than the average defined by its type. This variation can be driven by environmental conditions such as currents that affect different ocean areas with different magnitudes. The influence of these conditions on the battery usage is modelled with an energy factor $\energyFactor$ that is assigned to each research station $\station$ ($\energyFactor[\station]$) and route segment $\route{o}{d}$ ($\energyFactor[\route{o}{d}]$).

\subsection{Ship integration}
The model can also include a pre-existing ship ($\ship$) itinerary. The latter is made of different mooring or stationary locations, each of which is characterised by a time window $\timeWindow{\startTime{}}{\eendTime{}}$ when the ship arrives and leaves that location. Each MAV can be assigned to the ship via its waypoints for one of these situations:
\begin{itemize}
    \item \textbf{Have a fully recharged battery}: This happens in the scenario when the new battery will lead to better coverage of the remaining science missions. The minimum duration spent at the ship $\ship$ is $\minimum{\duration{\ship}}$;
    \item \textbf{Use the ship for transiting}: Because the ship is generally faster and safer in challenging environmental conditions. MAVs can transit via the ship to get to their destination faster when the proposed route would not allow them to reach their destinations within the specified time window or would lead to a battery resource constraint.
\end{itemize}
 
\subsection{Routes}
The stations $\setStation$, bases $\setBase$, pickup locations $\setPickup$ and the ship $\ship$ are interconnected via a set of routes $\setRoute$. Each route $\route{o}{d}$ is directed and has an origin \(o\) and a destination \(d\). In addition to $\energyFactor[\route{o}{d}]$, $\duration{\route{o}{d}}$ is the travel duration between \(o\) and \(d\).

As stated earlier, the ship is located at different positions across different time windows, making the duration and energy factors between the ship and any other location time-dependent. In addition to the origin \(o\) and destination \(d\) (where one of them is a ship), ship routes ($\route{\ship_{\timeWindow{}{}}}{d}$ or $\route{o}{\ship_{\timeWindow{}{}}}$) are also characterised by a time window $\timeWindow{}{}$ when the route is available.

\section{Summary of notations}
To transition seamlessly from the conceptual problem formulation to the rigorous mixed-integer linear programming model, we compile a comprehensive summary of all mathematical notations used throughout the rest of this paper. This reference is organised into three primary modules.

\subsubsection*{AUV types:}
\begin{itemize}
    \item $\setAUVType{}$: set of all MAV types;
    \item $\AUVType$: a single MAV type in $\setAUVType{}$;
    \item  $\setAUVType[\station]$: set of MAV types required by station $\station$
\end{itemize}

\subsection*{AUVs:}
\begin{itemize}
    \item $\setAUV$: the set of all available MAVs, of any type (the fleet);
    \item $\setAUV[\AUVType]$: the set of MAVs of type $\AUVType$;
    \item $\AUV$: a single MAV in $\setAUV$ (or $\setAUV[\AUVType]$);
\end{itemize}

\subsection*{Time Windows}
\begin{itemize}
    \item $\setTimeWindow$: set of all time windows;
    \item $\timeWindow{s}{e}$: a single time window, starting at \(s\) and ending at \(e\); 
    \item $\startTime{\timeWindow{}{}}$: start time of $\timeWindow{}{}$;
    \item $\eendTime{\timeWindow{}{}}$: end time of $\timeWindow{}{}$;
    \item $\timeWindow{}{}_a \leq \timeWindow{}{}_b \equiv \startTime{\timeWindow{}{}_a} \leq \startTime{\timeWindow{}{}_b}$
\end{itemize}

\subsection*{Locations:}
\begin{itemize}
    \item $\setStation$: set of all research stations;
    \item $\setStation[\AUVType]$: set of research stations to which the MAVs of type $\AUVType$ can be assigned;
    \item $\station$: a single station in $\setStation$;
    \item $\setBase$: the set of all ports of mobilisation for $\setAUV$;
    \item $\base{\AUV}$: the port of mobilisation for $\AUV$;
    \item $\setPickup$: the set of all ports of demobilisation for $\setAUV$;
    \item $\pickup{\AUV}$: the port of demobilisation for $\AUV$;
    \item $\ship[\timeWindow{}{}]$: the ship location within the time window \timeWindow{}{}
    \item $\ship[\setTimeWindow]$: the set of ship locations in each time window in $\setTimeWindow$
\end{itemize}

\subsection*{Routes:}
\begin{itemize}
    \item $\setRoute$: the set of all routes available for MAV movements;
    \item $\route{o}{d}$: a single static route in $\setRoute$ connecting an origin \(o\) and destination \(d\) points. They can be either $\station \in \setStation$, $\base{} \in \setBase$, $\pickup{} \in \setPickup$ or $\ship[\timeWindow{}{} \in \setTimeWindow]$;
    \item $\setRoute_{\ship}$: the set of ship routes, where the ship $\ship$ is one of the route endpoints.
\end{itemize}

\subsubsection*{Task Requirements:}
\begin{itemize}
    \item $\setTaskReq{}{}$: the set of all task requirements to be completed;
    \item $\setTaskReq{\station}{\AUVType}$: set of task requirements to be completed in $\station$ by MAVs of type $\AUVType$;
    \item $\taskReq{\station}{\AUVType}$: a single task in $\station$ requiring the type $\AUVType$;
    \item $\taskReqWeight{\taskReq{}{}}$: weight of the task requirement $\taskReq{}{}$
\end{itemize}

\subsection{Model Input Parameters}
\subsubsection*{AUV types:}
\begin{itemize}
    \item $\activeBatteryCons{\AUVType}$: average of battery consumption per time unit when the MAV of type $\AUVType$ is active;
    \item $\idleBatteryCons{\AUVType}$: average of battery consumption per time unit when the MAV of type $\AUVType$ is idle;
    \item $\batteryLevelThreshold{\AUVType}$: battery level threshold for each MAV of type $\AUVType$;
\end{itemize}

\subsubsection*{AUVs:}
\begin{itemize}
    \item $\initBatteryLevel{\AUV}$: the initial battery level for MAV $\AUV$;
\end{itemize}

\subsubsection*{Locations:}
\begin{itemize}
    \item $\energyFactor[\station]$: energy factor for each station $\station$;
\end{itemize}

\subsubsection*{Routes:}
\begin{itemize}
    \item $\duration{\route{o}{d}}$: travel duration between \(o\) and \(d\);
    \item $\energyFactor[\route{o}{d}]$: energy factor for each route $\route{o}{d}$;
\end{itemize}

\subsubsection*{Task Requirements:}
\begin{itemize}
    \item $\setAUVType[\station]$: the set of MAV types (a subset of $\setAUVType{}$) needed in the station $\station$; 
    \item $\duration{\taskReq{\station}{\AUVType}}$: duration of the task requirement $\taskReq{\station}{\AUVType}$ at the station $\station$ using MAVs of type $\AUVType$ ;
    \item $\requirement{\taskReq{\station}{\AUVType}}$: number of required MAVs from the type $\AUVType$ in the $\station$;
    \item $\earliestTime{\taskReq{}{}}$: earliest start time of a task requirement at $\station$ using MAVs of type $\AUVType$ ;
    \item $\latestTime{\taskReq{}{}}$: latest end time of a task requirement at $\station$ using MAVs of type $\AUVType$ ;
\end{itemize}

\subsubsection*{Departure \& Arrival:}
\begin{itemize}
    \item $\earliestDeployment{\base{}}$: earliest time for releasing the MAVs at $\base{}$;
    \item $\latestDeployment{\base{}}$: latest time for releasing the MAVs at $\base{}$;
    \item $\earliestRecovery{\pickup{}}$: earliest time for the MAVs to arrive at $\pickup{}$;
    \item $\latestRecovery{\pickup{}}$: latest time for the MAVs to arrive at $\pickup{}$;
\end{itemize}

\subsubsection*{Batteries \& Charging Facilities:}
\begin{itemize}
    \item $\maximum{\batteryLevel{}}:$ maximum of a battery level ($=100$, considering it is a percentage)
    \item $\minimum{\duration{\ship}}$: minimum duration spent in the ship $\ship$
\end{itemize}

\subsubsection*{Others:}
\begin{itemize}
    \item $\setTimeWindow$: set of all time windows;
    \item $\largeConstant$: a large positive constant
\end{itemize}

\subsection{Decision Variables}
\subsubsection*{Assignments:}
\begin{itemize}
    \item $\dAUVStation{\AUV}{\taskReq{}{}}$: binary variable for the assignment of the MAV $\AUV$ to the task $\taskReq{}{}$. \\ $\dAUVStation{\AUV}{\taskReq{}{}}=1$ if $\AUV$ is assigned to the task $\taskReq{}{}$, 0 otherwise;
    \item $\dAUVRoute{\AUV}{\route{}{}}$: binary variable representing the use of the route $\route{}{} \in \setRoute$ by the MAV $\AUV$.
    \begin{itemize}
        \item $\dAUVRoute{\AUV}{\route{o}{d}}=1$ if $\AUV$ is going from \(o\) to \(d\), 0 otherwise;
        \item $\dAUVRoute{\AUV}{\route{\ship_{\timeWindow{}{}}}{d}}=1$ if $\AUV$ is leaving the ship $\ship$ to \(d\) within $\timeWindow{}{}$;
        \item $\dAUVRoute{\AUV}{\route{o}{\ship_{\timeWindow{}{}}}}=1$ if $\AUV$ is heading to the ship $\ship$ from \(d\) within $\timeWindow{}{}$;
    \end{itemize}
    \item $\dAUVChargingFromTo{\AUV}{\timeWindow{}{}_a}{\timeWindow{}{}_b}$: binary variable for the assignment of the MAV $\AUV$ to the ship $\ship$, starting within the time window $\timeWindow{}{}_a$ and ending within $\timeWindow{}{}_b$.
    \\ $\dAUVChargingFromTo{\AUV}{\timeWindow{}{}_a}{\timeWindow{}{}_b}=1$ means $\AUV$ is staying in the ship between $\timeWindow{}{}_a$ and $\timeWindow{}{}_b$;
    \item $\dAUVChargingIn{\AUV}{\timeWindow{}{}}$: binary variable indicating whether the MAV $\AUV$ is heading to $\ship$ within $\timeWindow{}{}$ or not;
    \item $\dAUVChargingOut{\AUV}{\timeWindow{}{}}$: binary variable indicating whether the MAV $\AUV$ is leaving $\ship$ within $\timeWindow{}{}$ or not;
\end{itemize}

\subsubsection*{Timing:}
\begin{itemize}
    \item $\dAUVStationST{\AUV}{\taskReq{}{}}$: start time of the task $\taskReq{}{}$ by the MAV $\AUV$;
    \item $\dAUVRouteAT{\AUV}{\route{}{}}$: arrival time of MAV $\AUV$ at route ($\route{}{}$) destination;
\end{itemize}

\subsubsection*{Battery Levels:}
\begin{itemize}
    \item $\dBatteryLevelTask{\AUV}{\taskReq{\station}{\AUVType}}$: the battery level of MAV $\AUV$ at the end of task requirement $\taskReq{\station}{\AUVType}$
    \item $\dBatteryLevelRoute{\AUV}{\route{}{}}$: the battery level of MAV $\AUV$ at the end of the route $\route{}{}$
\end{itemize}

\section{Mathematical Model} 
\label{sec:model}

\subsection{Objective function}


The problem has five objectives defined in the following order:
\begin{enumerate}
    \item Maximise the number of MAV assignments to task requirements in the different stations:
    \\This is the main objective of the problem: maximising the number of science reviews comes before any other objective. A task weight $\taskReqWeight{\taskReq{}{}}$ is included to offer the possibility to prioritise some science missions on top of others when a full coverage cannot be attained due to operational constraints such as fleet size, time windows or MAV requirements;
    \begin{equation} \label{eq:obj_science}
        \max ~~ f_{\text{science}} = \sum_{\station \in \setStation} \sum_{\AUVType \in \setAUVType{\station}} \sum_{\AUV \in \setAUV[\AUVType]} \dAUVStation{\AUV}{\taskReq{\station}{\AUVType}} \times \frac{\taskReqWeight{\taskReq{}{}} }{\requirement{\taskReq{}{}}}
    \end{equation}
    
    \item Minimise the total number of MAV deployments:
    \\The model is not constrained to use the entire available MAV fleet. For cost reduction purposes, a smaller number of MAVs can be assigned if they can cover the same number of task requirements;
    \begin{equation} \label{eq:obj_deployment}
        \min ~~ f_{\text{deployment}} = \sum_{\AUVType \in \setAUVType} \sum_{\AUV \in \setAUV[\AUVType]} \sum_{d \in \setStation \cup \ship[\setTimeWindow]} \dAUVRoute{\AUV}{\route{\base{\AUV}}{d}}
    \end{equation}
    
    \item Minimise the number of MAV assignments to the ship:
    \\Each MAV assignment to the ship induces a penalty to the objective value. This will restrict the ship assignment to MAVs only when an equivalent or better solution cannot be found;
    \begin{equation} \label{eq:obj_docking}
        \min ~~ f_{\text{docking}} = \sum_{\AUVType \in \setAUVType} \sum_{\AUV \in \setAUV[\AUVType]} \times \sum_{\timeWindow{}{} \in \setTimeWindow} \dAUVChargingIn{\AUV}{\timeWindow{}{}}
    \end{equation}
    
    \item Minimise the battery usage for all the MAVs:
    \\This translates in the model by maximising the battery level for each deployed vehicle. We consider in the objective function the battery levels of MAVs when arriving at the ship or their respective pickup location;
    \begin{equation} \label{eq:obj_battery}
        \max ~~ f_{\text{battery}} = \sum_{\AUVType \in \setAUVType} \sum_{\AUV \in \setAUV[\AUVType]} \times \sum_{o \in \setStation \cup \{\base{\AUV}\}} \sum_{d \in \ship[\setTimeWindow] \cup \{\pickup{\AUV}\}} \dBatteryLevelRoute{\AUV}{\route{o}{d}}
    \end{equation}
    
    \item Minimise the arrival time to the pickup location (port of demobilisation): 
    \\As a final objective, the model minimises each deployed MAV's timespan by assigning it to its pickup location as soon as possible, within the specified time window.
    \begin{equation} \label{eq:obj_timespan}
        \min ~~ f_{\text{timespan}} = \sum_{\AUVType \in \setAUVType} \sum_{\AUV \in \setAUV[\AUVType]} \times \sum_{o \in \setStation \cup \ship[\setTimeWindow]} \dAUVRouteAT{\AUV}{\route{o}{\pickup{\AUV}}}
    \end{equation}
\end{enumerate}

The whole objective function combines the five above ones as in (\ref{eq:set_objective}). 
$\alpha, \beta, \gamma, \theta$ and $\omega$ are weight parameters to define the priority between each objective. In our paper, $\alpha = \largeConstant \times id_{\AUV} \times id_{\taskReq{}{}}$. We generate IDs for each vehicle and task requirement and include them in the objective function to avoid symmetries. Then $\beta = \largeConstant$, $\gamma = \maximum{\batteryLevel{}}^{2}$, $\theta = 10$ and $\omega = 1$.

\begin{equation} \label{eq:set_objective}
    \max ~~  \alpha \times f_{\text{science}} - ( \beta \times f_{\text{deployment}} + \gamma \times f_{\text{docking}} - \theta \times f_{\text{battery}} + \omega \times f_{\text{timespan}})
\end{equation}

\subsection{Flow constraints}

\begin{flalign} \label{eq:set_constraint_leaving}
        s.t. ~~  \sum_{d \in \setStation \cup \ship[\setTimeWindow]} \dAUVRoute{\AUV}{\route{\base{\AUV}}{d}} \leq 1
        ~,~ \forall \AUV \in \setAUV  &&
\end{flalign}

\begin{flalign} \label{eq:set_constraint_pickup}
        \sum_{o \in \setStation \cup \ship[\setTimeWindow]} \dAUVRoute{\AUV}{\route{o}{\pickup{\AUV}}} \leq 1
        ~,~ \forall \AUV \in \setAUV  &&
\end{flalign}

\begin{flalign} \label{eq:set_constraint_assignment}
        \sum_{\AUV \in \setAUV[\AUVType]} \dAUVStation{\AUV}{\taskReq{\station}{\AUVType}} \leq \requirement{\taskReq{\station}{\AUVType}}
        ~,~ \forall \station \in \setStation ~ \forall \AUVType \in \setAUVType[\station] ~ \mid \exists \taskReq{\station}{\AUVType} \in \setTaskReq{}{} &&
\end{flalign}

\begin{flalign} {\label{eq:set_constraint_intratask_incoming}}
    \begin{aligned}
            \dAUVStation{\AUV}{\taskReq{\station}{\AUVType}} 
            & = \sum_{o \in (\setStation \setminus \{\station\} )\cup \ship[\setTimeWindow] \cup \{\base{\AUV}\}} \dAUVRoute{\AUV}{\route{o}{\station}} ~,~ &\\
            & \forall \station \in \setStation  ~,~ \forall \AUVType \in \setAUVType[\station] ~,~ \forall \AUV \in \setAUV[\AUVType]
    \end{aligned}
    &&
\end{flalign}

\begin{flalign} {\label{eq:set_constraint_intratask_outgoing}}
    \begin{aligned}
            \dAUVStation{\AUV}{\taskReq{\station}{\AUVType}}  
            & = \sum_{d \in (\setStation \setminus \{\station\}) \cup \ship[\setTimeWindow] \cup \{\pickup{\AUV}\}} \dAUVRoute{\AUV}{\route{\station}{d}} ~,~ &\\
            & \forall \station \in \setStation  ~,~ \forall \AUVType \in \setAUVType[\station] ~,~ \forall \AUV \in \setAUV[\AUVType]
    \end{aligned}
    &&
\end{flalign}

\begin{flalign} \label{eq:set_constraint_intratask_charging_incoming}
     \dAUVChargingIn{\AUV}{\timeWindow{}{}} = \sum_{o \in \setStation \cup \{\base{\AUV}\}} \dAUVRoute{\AUV}{\route{o}{\ship_{\timeWindow{}{}}}} ~,~ \forall \AUV \in \setAUV ~,~ \forall \timeWindow{}{} \in \setTimeWindow && 
\end{flalign}

\begin{flalign} \label{eq:set_constraint_intratask_charging_outgoing}
     \dAUVChargingOut{\AUV}{\timeWindow{}{}} = \sum_{d \in \setStation \cup \{\pickup{\AUV}\}} \dAUVRoute{\AUV}{\route{\ship_{\timeWindow{}{}}}{d}} ~,~ \forall \AUV \in \setAUV ~,~ \forall \timeWindow{}{} \in \setTimeWindow &&
\end{flalign}

\begin{flalign} \label{eq:set_constraint_assignment_cf.1}
        \dAUVChargingIn{\AUV}{\timeWindow{}{}} =  \sum_{\timeWindow{}{}_t \in \setTimeWindow : \timeWindow{}{}_t \geq \timeWindow{}{}} \dAUVChargingFromTo{\AUV}{\timeWindow{}{}}{\timeWindow{}{}_t} ~,~ \forall \AUV \in \setAUV{} ~,~ \forall \timeWindow{}{} \in \setTimeWindow &&
\end{flalign}

\begin{flalign} \label{eq:set_constraint_assignment_cf.2}
        \dAUVChargingOut{\AUV}{\timeWindow{}{}} =  \sum_{\timeWindow{}{}_t \in \setTimeWindow : \timeWindow{}{}_t \leq \timeWindow{}{}} \dAUVChargingFromTo{\AUV}{\timeWindow{}{}_t}{\timeWindow{}{}}  ~,~ \forall \AUV \in \setAUV{} ~,~ \forall \timeWindow{}{} \in \setTimeWindow &&
\end{flalign}

As mentioned earlier in the paper, we deploy each MAV only if it can lead to a new optimal solution. Constraints (\ref{eq:set_constraint_leaving}) state that each MAV will be leaving its base location via, at most, one route. The same holds when an MAV has to be recovered at its pickup location, as verified by (\ref{eq:set_constraint_pickup}).
Constraints (\ref{eq:set_constraint_assignment}) state that, for each task requirement $\taskReq{\station}{\AUVType}$, the number of assigned MAVs of type $\AUVType$ to the station $\station$ should not exceed the number of required MAVs by the task. 
When an MAV is assigned to a task, exactly one route should be activated to come to the station, and one to leave it. The two are controlled respectively by the constraints (\ref{eq:set_constraint_intratask_incoming}) and (\ref{eq:set_constraint_intratask_outgoing}), which also ensure that if an MAV is not assigned to a task/station, it will not traverse that location. 
Similarly, if an MAV is assigned to the ship within one of the available time windows, there should be only one incoming (resp. outgoing) route to the ship at that location (or time window) activated, which is verified by constraints (\ref{eq:set_constraint_intratask_charging_incoming}) (resp. (\ref{eq:set_constraint_intratask_charging_outgoing})).
Since the model allows for each MAV to be assigned multiple times, constraints (\ref{eq:set_constraint_assignment_cf.1}) ensure that, if the MAV is assigned to the ship within a time window $\timeWindow{}{}$, then it would have to leave the ship within the same time window or one of the following ones. The same goes for constraints (\ref{eq:set_constraint_assignment_cf.2}), making sure that if an MAV is leaving the ship within a time window, then it should have been assigned to the ship within the same one or earlier.

\subsection{Time constraints}

\begin{flalign} {\label{eq:set_constraint_arrival_time_from_base}} 
    \begin{aligned}
            (\earliestDeployment{\base{\AUV}} + \duration{\route{\base{\AUV}}{d}}) \times \dAUVRoute{\AUV}{\route{\base{\AUV}}{d}} & \leq \dAUVRouteAT{\AUV}{\route{\base{\AUV}}{d}} \leq (\latestDeployment{\base{\AUV}} + \duration{\route{\base{\AUV}}{d}}) \times \dAUVRoute{\AUV}{\route{\base{\AUV}}{d}} ~,~ &\\
            & \forall d \in \setStation \cup \ship[\setTimeWindow] ~,~ \forall \AUV \in \setAUV
    \end{aligned}
    &&
\end{flalign}
\begin{flalign} \label{eq:set_constraint_arrival_time_rdv_bounds} 
\begin{aligned}
    \earliestRecovery \times \dAUVRoute{\AUV}{\route{o}{\pickup{\AUV}}} \leq \dAUVRouteAT{\AUV}{\route{o}{\pickup{\AUV}}} \leq & \latestRecovery \times \dAUVRoute{\AUV}{\route{o}{\pickup{\AUV}}} ~,~ &\\
    & \forall o \in \setStation \cup \ship[\setTimeWindow] ~,~ \forall \AUVType \in \setAUVType[\station] ~,~ \forall \AUV \in \setAUV[\AUVType] 
\end{aligned}
&&
\end{flalign}


\begin{flalign} {\label{eq:set_constraint_arrival_time}} 
\begin{aligned} 
        \dAUVRouteAT{\AUV}{\route{\station_o}{d}} & \geq \dAUVStationST{\AUV}{\taskReq{\station_o}{\AUVType}} + 
        \duration{\taskReq{\station_o}{\AUVType}} + \duration{\route{\station_o}{d}} - K \times (1 - \dAUVRoute{\AUV}{\route{\station_o}{\station_d}}), & \\
        & \forall \station_o \in \setStation , \forall d \in \setStation \cup \{\pickup{\AUV}\} , \forall \AUVType \in \setAUVType[\station_o] ~,~ \forall \AUV \in \setAUV[\AUVType] 
\end{aligned} 
        &&
\end{flalign}

\begin{flalign} \label{eq:set_constraint_arrival_time_to_cf} 
\begin{aligned} 
        \dAUVRouteAT{\AUV}{\route{\station_o}{\ship_{\timeWindow{}{}}}} & \geq \dAUVStationST{\AUV}{\taskReq{\station_o}{\AUVType}} + 
        \duration{\taskReq{\station_o}{\AUVType}} + (\duration{\route{\station_o}{\ship_{\timeWindow{}{}}}} \times \dAUVChargingIn{\AUV}{\timeWindow{}{}}) - K \times (1 - \dAUVRoute{\AUV}{\route{\station_o}{\ship_{\timeWindow{}{}}}}), & \\
        & \forall \station_o \in \setStation  , \forall \AUVType \in \setAUVType{\station_o} ~,~ \forall \AUV \in \setAUV[\AUVType] ~,~ \forall \timeWindow{}{} \in \setTimeWindow
\end{aligned} 
        &&
\end{flalign}


\begin{flalign} \label{eq:set_constraint_arrival_time_from_charging}
\begin{aligned}
    \dAUVRouteAT{\AUV}{\route{\ship_{\timeWindow{}{}}}{d}} & \geq \sum_{o \in \setStation \cup \{\base{\AUV}\}} (\dAUVRouteAT{\AUV}{\route{o}{\ship_{\timeWindow{}{}}}} + \dAUVRoute{\AUV}{\route{o}{\ship_{\timeWindow{}{}}}} \times \minimum{\duration{\ship}})  &\\
    & + \duration{\route{\ship_{\timeWindow{}{}}}{d}} \times \dAUVRoute{\AUV}{\route{\ship_{\timeWindow{}{}}}{d}} - \largeConstant \times (1 - \dAUVRoute{\AUV}{\route{\ship_{\timeWindow{}{}}}{d}}) ~,~ &\\
    & \forall d \in \setStation \cup \{\pickup{\AUV}\} ~,~  \forall \AUV \in \setAUV ~,~ \forall \timeWindow{}{} \in \setTimeWindow
\end{aligned}
    &&
\end{flalign}
\begin{flalign} \label{eq:set_constraint_arrival_time_bounds_from_charging}
    \begin{aligned}
        ((\startTime{\timeWindow{}{}} + \duration{\route{\ship_{\timeWindow{}{}}}{d}}) \times \dAUVRoute{\AUV}{\route{\ship_{\timeWindow{}{}}}{d}}) & \leq \dAUVRouteAT{\AUV}{\route{\ship_{\timeWindow{}{}}}{d}} \leq  ((\eendTime{\timeWindow{}{}} + \duration{\route{\ship_{\timeWindow{}{}}}{d}}) \times \dAUVRoute{\AUV}{\route{\ship_{\timeWindow{}{}}}{d}}) ~,~ &\\
        & \forall d \in \setStation \cup \{\pickup{\AUV}\} ~,~ \forall \AUV \in \setAUV ~,~ \forall \timeWindow{}{} \in \setTimeWindow
    \end{aligned}
     &&
\end{flalign}

\begin{flalign} \label{eq:set_constraint_arrival_time_bounds_to_charging}
    \begin{aligned}
        \startTime{\timeWindow{}{}} \times \dAUVRoute{\AUV}{\route{\station}{\ship_{\timeWindow{}{}}}} & \leq \dAUVRouteAT{\AUV}{\route{\station}{\ship_{\timeWindow{}{}}}} \leq  \latestTime{\timeWindow{}{}} \times \dAUVRoute{\AUV}{\route{\station}{\ship_{\timeWindow{}{}}}} ~,~ &\\
        & \forall \station \in \setStation \cup \{\base{\AUV}\} ~,~ \forall \AUV \in \setAUV ~,~ \forall \timeWindow{}{} \in \setTimeWindow
    \end{aligned}
    &&
\end{flalign}

\begin{flalign} {\label{eq:set_constraint_arrival_time_ub.1}}
\begin{aligned}
    \dAUVRouteAT{\AUV}{\route{o}{d}} \leq & \dAUVRoute{\AUV}{\route{o}{d}} \times \largeConstant ~,~ \forall \AUV \in \setAUV ~,~ &\\ 
    & \forall (o, d) \in (\setStation \cup \{\base{\AUV}\}) \times (\setStation \cup \{\pickup{\AUV}\}) \mid o \neq d 
\end{aligned}
&&
\end{flalign}
\begin{flalign} {\label{eq:set_constraint_arrival_time_ub.2}}
\begin{aligned}
    \dAUVRouteAT{\AUV}{\route{o_{\timeWindow{}{}}}{d_{\timeWindow{}{}}}} \leq & \dAUVRoute{\AUV}{\route{o_{\timeWindow{}{}}}{d_{\timeWindow{}{}}}} \times \largeConstant ~,~ \forall \AUV \in \setAUV ~,~ \forall \timeWindow{}{} \in \setTimeWindow  ~,~ &\\ 
    & \forall (o, d) \in (\setStation \cup \{\base{\AUV}, \ship\}) \times (\setStation \cup \{\pickup{\AUV}, \ship\}) \mid o \neq d \wedge (o = \ship \vee d = \ship) 
\end{aligned}
&&
\end{flalign}
\begin{flalign} {\label{eq:set_constraint_start_time_task.1}}
\begin{aligned}  
        \dAUVStationST{\AUV}{\taskReq{\station_d}{\AUVType}} \geq & \sum_{o \in \setStation \cup \{\base{\AUV}\}} \dAUVRouteAT{\AUV}{\route{o}{\station_d}} + \sum_{\timeWindow{}{} \in \setTimeWindow} \dAUVRouteAT{\AUV}{\route{\ship_{\timeWindow{}{}}}{\station_d}} ~,~ &\\ 
        & \forall \station_d \in \setStation ~,~ \forall \AUVType \in \setAUVType[\station_d] ~,~ \forall \AUV \in \setAUV[\AUVType]
\end{aligned}
&&
\end{flalign}

\begin{flalign} {\label{eq:set_constraint_start_time_task.3}}
        \earliestTime{\taskReq{\station}{\AUVType}} \times \dAUVStation{\AUV}{\taskReq{\station}{\AUVType}} \leq \dAUVStationST{\AUV}{\station} \leq \latestTime{\taskReq{\station}{\AUVType}}-\duration{\taskReq{\station}{\AUVType}}
         ~,~ \forall \station \in \setStation ~,~ \forall \AUVType \in \setAUVType{\station} ~,~ \forall \AUV \in \setAUV[\AUVType]  &&
\end{flalign}


Constraints (\ref{eq:set_constraint_arrival_time_from_base}) ensure that each deployed MAV leaves the base between the earliest time and the latest time of deployment. Deployed MAVs must then arrive at their respective pickup locations between the earliest and latest recovery times. This is asserted by constraints (\ref{eq:set_constraint_arrival_time_rdv_bounds}).
Constraints (\ref{eq:set_constraint_arrival_time}, \ref{eq:set_constraint_arrival_time_to_cf}) calculate the arrival time when the route origin is a research station, with (\ref{eq:set_constraint_arrival_time_to_cf}) being for the case when the destination is the ship. If the latter is a starting point, and whatever the type of the destination, the arrival time for this route is calculated by the constraints (\ref{eq:set_constraint_arrival_time_from_charging}). 
The arrival time for every ship route is constrained by the time window when the ship is available at its respective location. Constraints  (\ref{eq:set_constraint_arrival_time_bounds_from_charging}, \ref{eq:set_constraint_arrival_time_bounds_to_charging}) set the upper and lower bounds for ship route arrival times, respectively, when leaving her or heading to.
The arrival time for each MAV at an unused route is set at 0, as defined by the constraints (\ref{eq:set_constraint_arrival_time_ub.1}, \ref{eq:set_constraint_arrival_time_ub.2}).
As stated by the constraints (\ref{eq:set_constraint_start_time_task.1}), the start time of a task by an MAV is equivalent to its arrival time at the station, and it should be within a time window to allow the task to be completed, as defined by the constraints (\ref{eq:set_constraint_start_time_task.3}).
\subsection{Battery levels}

\begin{flalign} {\label{eq:set_constraint_battery_level_task_bounds}}
\begin{aligned}
    \batteryLevelThreshold{\AUVType} \times \dAUVStation{\AUV}{\taskReq{\station}{\AUVType}} \leq & \dBatteryLevelTask{\AUV}{\taskReq{\station}{\AUVType}} \leq \maximum{\batteryLevel{}} \times  \dAUVStation{\AUV}{\taskReq{\station}{\AUVType}} ~,~&\\ 
    & \forall \station \in \setStation ~,~ \forall \AUVType \in \setAUVType{\station} ~,~ \forall \AUV \in \setAUV[\AUVType]
\end{aligned}
    &&
\end{flalign}
\begin{flalign} {\label{eq:set_constraint_battery_level_route_bounds.1}}
\begin{aligned}
    \batteryLevelThreshold{\AUVType} \times \dAUVRoute{\AUV}{\route{o}{d}} \leq & \dBatteryLevelRoute{\AUV}{\route{o}{d}} \leq \maximum{\batteryLevel{}} \times  \dAUVRoute{\AUV}{\route{o}{d}} ~,~&\\ 
    & \forall (o, d) \in (\setStation \cup \ship[\setTimeWindow] \cup \{\base{\AUV}\}) \times (\setStation \cup \ship[\setTimeWindow] \cup \{\pickup{\AUV}\}) \mid o \neq d ~,~ \forall \AUV \in \setAUV
\end{aligned}
    &&
\end{flalign}
\begin{flalign} {\label{eq:set_constraint_battery_level_task}}
    \begin{aligned}
        \dBatteryLevelTask{\AUV}{\taskReq{\station}{\AUVType}} \leq & \sum_{o \in \setStation[\AUVType] \cup \{\base{\AUV}\}} \dBatteryLevelRoute{\AUV}{\route{o}{\station}} + \sum_{\timeWindow{}{} \in \setTimeWindow} \dBatteryLevelRoute{\AUV}{\route{\ship_{\timeWindow{}{}}}{\station}} - \duration{\taskReq{\station}{\AUVType}} \times \activeBatteryCons{\AUVType} \times \energyFactor[\station] \times \dAUVStation{\AUV}{\taskReq{\station}{\AUVType}} &\\
        & - (\dAUVStationST{\AUV}{\taskReq{\station}{\AUVType}} - \sum_{\station_o \in \setStation[\AUVType]} \dAUVRouteAT{\AUV}{\route{o}{\station}} - \sum_{\timeWindow{}{} \in \setTimeWindow} \dAUVRouteAT{\AUV}{\route{\ship_{\timeWindow{}{}}}{\station}} ) \idleBatteryCons{\AUVType} \times \energyFactor[\station] ~, &\\
        & \forall \station \in \setStation ~,~ \forall \AUVType \in \setAUVType[\station] ~,~ \forall \AUV \in \setAUV[\AUVType]
    \end{aligned}
    &&
\end{flalign}

\begin{flalign} \label{eq:set_constraint_battery_level_route1.1}
\begin{aligned}
    \dBatteryLevelRoute{\AUV}{\route{\station_o}{d}} \leq & \dBatteryLevelTask{\AUV}{\taskReq{\station_o}{\AUVType}} - \duration{\route{\station_o}{d}} \times \idleBatteryCons{\AUVType} \times \energyFactor[\route{\station_o}{d}] \times \dAUVRoute{\AUV}{\route{\station_o}{d}} &\\
    & - ((\dAUVRouteAT{\AUV}{\route{\station_o}{d}} - \duration{\route{\station_o}{d}}) - (\dAUVStationST{\AUV}{\taskReq{\station_o}{\AUVType}} + \duration{\taskReq{\station_o}{\AUVType}})) \times \idleBatteryCons{\AUVType} \times \energyFactor[\station_o] &\\
    & \forall \station_o \in \setStation ~,~ \forall d \in (\setStation \cup \ship[\setTimeWindow] \cup \{\pickup{\AUV}\}) ~,~ \forall \AUVType \in \setAUVType[\station] ~,~ \forall \AUV \in \setAUV[\AUVType]
\end{aligned}
&&
\end{flalign}
\begin{flalign} \label{eq:set_constraint_battery_level_route2}
\begin{aligned}
    \dBatteryLevelRoute{\AUV}{\route{\ship_{\timeWindow{}{}}}{\station_d}} \leq & \maximum{\batteryLevel{\AUV}} - \durationCF{\route{\ship_{\timeWindow{}{}}}{\station_d}}{\timeWindow{}{}} \times \energyFactorCF{\route{\ship_{\timeWindow{}{}}}{\station_d}}{\timeWindow{}{}} \times \idleBatteryCons{\AUVType} \times \dAUVRoute{\AUV}{\route{\ship_{\timeWindow{}{}}}{\station_d}} &\\
    & \forall \station_d \in \setStation \cup \{\pickup{\AUV}\} ~,~ \forall \AUVType \in \setAUVType{\station} ~,~ \forall \AUV \in \setAUV[\AUVType] ~,~  \forall \timeWindow{}{} \in \setTimeWindow
\end{aligned}
&&
\end{flalign}
\begin{flalign} \label{eq:set_constraint_battery_level_route3.1}
\begin{aligned}
    \dBatteryLevelRoute{\AUV}{\route{\base{\AUV}}{\station_d}} \leq & \initBatteryLevel{\AUV} - \duration{\route{\base{\AUV}}{\station_d}} \times \idleBatteryCons{\AUVType} \times \energyFactor[\route{\base{\AUV}}{\station_d}] \times \dAUVRoute{\AUV}{\route{\base{\AUV}}{\station_d}} ~,~ &\\
    & \forall d \in (\setStation \cup \ship_{\setTimeWindow}) ~,~ \forall \AUVType \in \setAUVType[\station] ~,~ \forall \AUV \in \setAUV[\AUVType]
\end{aligned}
&&
\end{flalign}

Constraints (\ref{eq:set_constraint_battery_level_task_bounds}) (resp. (\ref{eq:set_constraint_battery_level_route_bounds.1})) define the battery level bounds when at the end of any task (resp. route).
The battery level at the end of a task is calculated by the constraints (\ref{eq:set_constraint_battery_level_task}).
The calculation of the battery level at the end of the routes is based on the following three cases: (i) when leaving a research station (\ref{eq:set_constraint_battery_level_route1.1}), (ii) when leaving the ship (\ref{eq:set_constraint_battery_level_route2}) and (iii) when leaving the base, meaning the beginning of the deployment (\ref{eq:set_constraint_battery_level_route3.1}).
\subsection{Variable bound definitions}

\begin{flalign} {\label{eq:set_variables_dAUVStation}}
        \dAUVStation{\AUV}{\taskReq{\station}{\AUVType}} \in \{0,1\}
        ~,~ \forall \station \in \setStation ~,~ \forall \AUVType \in \setAUVType{\station} ~,~ \forall \AUV \in \setAUV[\AUVType]  &&
\end{flalign}

\begin{flalign} {\label{eq:set_variables_dAUVRoute}} 
        \dAUVRoute{\AUV}{\route{o}{d}} \in \{0,1\}
        ~,~ \forall \AUV \in \setAUV ~,~ \forall \route{o}{d} \in \setRoute &&
\end{flalign}

\begin{flalign} {\label{eq:set_variables_dAUVChargingFromTo}} 
        \dAUVChargingFromTo{\AUV}{\timeWindow{}{}_i}{\timeWindow{}{}_j} \in \{0,1\}
        ~,~ \forall \AUV \in \setAUV ~,~ \forall (\timeWindow{}{}_i,\timeWindow{}{}_i) \in \setTimeWindow^2 &&
\end{flalign}

\begin{flalign} {\label{eq:set_variables_dAUVChargingInOut}}
        \dAUVChargingIn{\AUV}{\timeWindow{}{}}, \dAUVChargingOut{\AUV}{\timeWindow{}{}} \in \{0,1\}
        ~,~ \forall \AUV \in \setAUV ~,~ \forall \timeWindow{}{} \in \setTimeWindow &&
\end{flalign}

\begin{flalign} {\label{eq:set_variables_dBatteryLevelRoute}}
        0 \leq \dBatteryLevelRoute{\AUV}{\route{o}{d}} \leq \maximum{\batteryLevel{}}
         ~,~ \forall \AUV \in \setAUV ~,~ \forall \route{o}{d} \in \setRoute &&
\end{flalign}

\begin{flalign} {\label{eq:set_variables_dBatteryLevelRoute}}
        0 \leq \dBatteryLevelTask{\AUV}{\taskReq{\station}{\AUVType}} \leq \maximum{\batteryLevel{}}
         ~,~ \forall \AUVType \in \setAUVType ~,~ \forall \AUV \in \setAUV[\AUVType] ~,~ \forall \station \in \setStation[\AUVType] &&
\end{flalign}

\begin{flalign} {\label{eq:set_variables_dAUVRouteAT}}
        0 \leq \dAUVRouteAT{\AUV}{\route{o}{d}} \leq \latestTime{\pickup{\AUV}}
         ~,~ \forall \AUV \in \setAUV ~,~ \forall \route{o}{d} \in \setRoute \mid o \neq \ship \wedge d \neq \ship  &&
\end{flalign}

\begin{flalign} {\label{eq:set_variables_dAUVRouteAT2}}
        0 \leq \dAUVRouteAT{\AUV}{\route{\ship_{\timeWindow{}{}}}{d}} \leq \eendTime{\timeWindow{}{}} + \duration{\route{\ship_{\timeWindow{}{}}}{d}}
         ~,~ \forall \AUV \in \setAUV ~,~ \forall \route{\ship_{\timeWindow{}{}}}{d} \in \setRoute_{\ship}  &&
\end{flalign}

\begin{flalign} {\label{eq:set_variables_dAUVRouteAT3}}
        0 \leq \dAUVRouteAT{\AUV}{\route{o}{\ship_{\timeWindow{}{}}}} \leq \eendTime{\timeWindow{}{}}
         ~,~ \forall \AUV \in \setAUV ~,~ \forall \route{o}{\ship_{\timeWindow{}{}}} \in \setRoute_{\ship}  &&
\end{flalign}

\begin{flalign} {\label{eq:set_variables_dAUVStationST}}
        0 \leq \dAUVStationST{\AUV}{\taskReq{}{}} \leq \latestTime{\taskReq{}{}} - \duration{\taskReq{}{}}
         ~,~ \forall \AUV \in \setAUV ~,~ \forall \taskReq{}{} \in \setTaskReq{}{}  &&
\end{flalign}

\begin{flalign} {\label{eq:set_variables_dAUVStationST2}}
        0 \leq \dAUVStationST{\AUV}{\ship[\timeWindow{}{}]} \leq \eendTime{\timeWindow{}{}}
         ~,~ \forall \AUV \in \setAUV ~,~ \forall \timeWindow{}{} \in \setTimeWindow  &&
\end{flalign}


\section{Case studies and Results}
\label{sec:case_studies}


In this section, we present four case studies that leverage the proposed MILP model to plan MAV fleets, maximising the number of science tasks completed while minimising the number of assets deployed and their battery consumption. The case studies vary in location, science requirements, and fleet size; each is inspired by real MAV science missions but deliberately expanded with additional task requirements and fleet assets. 
This scaling is intentional: as fleet size and the density of scientific commitments grow, manually planning a feasible deployment plan quickly becomes impractical, and even experienced human planners struggle to find good, let alone optimal, solutions within a reasonable time budget. 
The first case study establishes a baseline, shoreside-only deployment; the second introduces a hybrid deployment approach combining shoreside locations with a research vessel operating following an existing itinerary, to test the model’s handling of ship-assisted battery swapping and transit; the third evaluates the model at a near-national scale, coordinating a shared fleet of 100 vehicles across three different organisations; and the fourth demonstrates the model’s applicability across a diverse, worldwide portfolio of seven independently-planned deployments spanning four continents.

 All four case studies are solved using SCIP \citep{SCIPOptSuite10}, a non-commercial solver for mixed-integer programming (MIP) and mixed-integer nonlinear programming (MINLP).

\subsection{Inputs and outputs}

\subsubsection{Data sources}
The case studies presented and solved in this paper are inspired by real-world science missions taken from the Marine Facilities Planning (MFP) platform \footnote{Maas Software Engineering. Marine Facilities Planning (MFP), 2026, https://www.marinefacilitiesplanning.com/.}.
The MFP provides a record of Ship time and Marine Equipment (SME) applications, as well as Autonomous Deployment (ADF) forms, which we use as a basis for our case studies.

\textbf{Task requirements}:
ADF forms contain the set of task requirements, the locations of mobilisation (bases) and demobilisation (pickup location or rendezvous points) ports, and the research stations. The form specifies the MAV type, the number of vehicles required, and the planned deployment and recovery dates. 

\textbf{Ship itinerary}:
On the other hand, SMEs are used to build a ship itinerary that matches the task requirements to allow the MAVs to be assigned to the ship to cover additional task requirements. As it is a model with a single ship, these use cases focus on SMEs operated by the Royal Research Ship (RRS) Sir David Attenborough (SDA).

\textbf{AUVs and their types}:
This model considers the use and deployment of a fleet composed of Slocum Glider and Autosub Long Range 1500 (ALR1500) MAVs. The choice of these two is due to the availability of their characteristics within the organisation and other partners, which we can rely on for better route cost estimations.

In this model, we consider an MAV type to be a combination of either the Slocum Glider or the ALR1500 with different sensor configurations.

\textbf{Route costs}:
While the ADF forms provide the necessary locations to solve each case study, they don't include route costs between origins and destinations. PolarRoute \citep{smith2022autonomous}, an advanced AI- and data-driven maritime route-planning software, is used to generate travel time and battery consumption for the entire set of routes, including ship routes. We generate costs for the Slocum Glider and ALR1500, and for generating ship itineraries in some case studies.  

\subsubsection{Model outputs}

We generate different outputs for each case study in this work. They are either textual or visualisations:
\begin{itemize}
    \item \textbf{The case study plan}, a CSV file format with all the activities for each deployed MAV, starting from the deployment (departure) till the recovery (arrival). Activities can be science reviews, passages or ship assignments (battery swapping or transits). Each activity has a start time, an end time, and battery levels at the start and end.
    \item \textbf{Completion rates} of the task requirements. As explained earlier in this paper, each science task requires a number of MAVs of a specific type/configuration. The completion rate indicates how much is covered by the proposed plan and how much remains to be covered. Each rate maps a research station with an AVU type.
    \item A list of \textbf{unused MAVs}, which can be safely removed from the fleet before the start of the operations.
    \item A \textbf{timeline} of MAV activities generated based on the CSV plan. Each row of the timeline (example Figure \ref{fig:timeline_case_study1}) shows all the activities for a single MAV. The battery level variation is added on top of activities that consume battery.
    \item \textbf{Histogram} of science data collected by all the deployed MAVs. This data is grouped by instrument or sensor and counted using Science Days as the unit of measure.
    \item \textbf{Heatmap} to visualise the completion rates. It maps MAV types to research stations.
    \item \textbf{Maps} with all the different locations and the ship itinerary. Each map is a snapshot of the proposed plan at a specific date.
\end{itemize}

While we generate all the above outputs for each case study, we present the relevant visualisations below. All the visualisations are generated using the plotly package \citep{Kruchten_An_interactive_open-source_2026}, are interactive, and provide more information than the screenshots included in this paper.
\subsection{Case studies}

\subsection{Deployment around Rothera, Antarctica} \label{subsubsec:case_study1} \label{sssec:rothera_results}

\begin{figure}[h]
    \centering
    \includegraphics[width=0.9\linewidth]{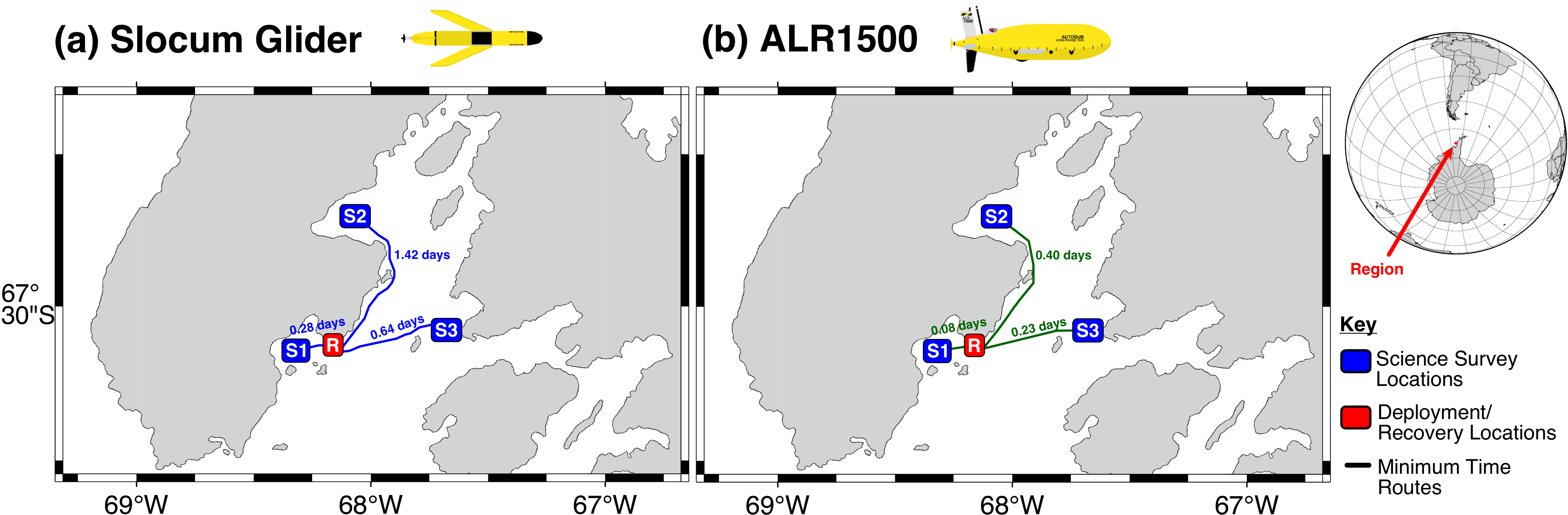}
    \caption{Shoreside deployment from Rothera, Antarctica. (a) Represents example routes from the deployment/recovery location Rothera (R) to each of the research stations (S1, S2 and S3) for Slocum marine gliders. The transit time in days for each route is labelled. (b) Represents the same as panel (a) but for the AL1500.}
    \label{fig:case_study1}
\end{figure}

This first case study is inspired by science surveys conducted by the British Antarctic Survey (BAS) around Rothera, the organisation's largest Antarctic research station and a key hub for Southern Ocean monitoring. It demonstrates the deployment and recovery of a mixed fleet of MAVs from a single location, with vehicles tasked to science surveys across three different research stations (S1, S2 and S3). The autonomous fleet consists of two marine platforms: Slocum Gliders (10 in total) and ALR1500 vehicles (2 in total), the latter moving at a faster speed than the Slocum Gliders (routes shown in Figure \ref{fig:case_study1}), with further differences between assets arising from their sensor instrumentation (Table \ref{tab:fleet_case_study1}). For each research station, the model is provided with the science surveys that need to be completed, along with the sensor requirements, the number of MAVs needed, the task duration, and the earliest-start/latest-end timing restrictions (Table \ref{tab:tasks_case_study1}); in total, 11 task requirements are distributed across the 3 research stations and 4 distinct MAV types, creating a combinatorially rich assignment problem despite the comparatively short distances involved.

\begin{table}[h]
\resizebox{\textwidth}{!}{%
\begin{tabular}{|cc|c|c|c|}
\hline
\multicolumn{2}{|c|}{\textbf{AUV Types}} & \multirow{2}{*}{\textbf{Base Location}} & \multirow{2}{*}{\textbf{Pickup Location}} & \multirow{2}{*}{\textbf{\begin{tabular}[c]{@{}c@{}}Initial Battery\\ Level (\%)\end{tabular}}} \\ \cline{1-2}
\multicolumn{1}{|c|}{\textbf{Autonomous Platform}} & \textbf{Sensors} &  &  &  \\ \hline
\multicolumn{1}{|c|}{Slocum} & CTD+microstructure+ADCP & Rothera & Rothera & 80 \\ \hline
\multicolumn{1}{|c|}{Slocum} & CTD+microstructure+ADCP & Rothera & Rothera & 60 \\ \hline
\multicolumn{1}{|c|}{Slocum} & CTD+microstructure+ADCP & Rothera & Rothera & 100 \\ \hline
\multicolumn{1}{|c|}{Slocum} & CTD+microstructure+ADCP & Rothera & Rothera & 88 \\ \hline
\multicolumn{1}{|c|}{Slocum} & CTD+ADCP & Rothera & Rothera & 80 \\ \hline
\multicolumn{1}{|c|}{Slocum} & CTD+ADCP & Rothera & Rothera & 100 \\ \hline
\multicolumn{1}{|c|}{Slocum} & CTD+ADCP & Rothera & Rothera & 100 \\ \hline
\multicolumn{1}{|c|}{Slocum} & CTD+optode+PAR+Ecopuck & Rothera & Rothera & 100 \\ \hline
\multicolumn{1}{|c|}{Slocum} & CTD+optode+PAR+Ecopuck & Rothera & Rothera & 100 \\ \hline
\multicolumn{1}{|c|}{Slocum} & CTD+optode+PAR+Ecopuck & Rothera & Rothera & 100 \\ \hline
\multicolumn{1}{|c|}{ALR1500} & UVP6+EK80+EcoBBRTD & Rothera & Rothera & 100 \\ \hline
\multicolumn{1}{|c|}{ALR1500} & UVP6+EK80+EcoBBRTD & Rothera & Rothera & 100 \\ \hline
\end{tabular}%
}
\caption{Fleet details for deployment around Rothera}
\label{tab:fleet_case_study1}
\end{table}

\begin{table}[h]
\resizebox{\textwidth}{!}{%
\begin{tabular}{|c|c|c|c|c|c|}
\hline
\textbf{\begin{tabular}[c]{@{}c@{}}Research \\ Stations\end{tabular}} & \textbf{AUV Types} & \textbf{\begin{tabular}[c]{@{}c@{}}Number \\ Required\end{tabular}} & \textbf{\begin{tabular}[c]{@{}c@{}}Task Duration \\ (days)\end{tabular}} & \textbf{Earliest Start} & \textbf{Latest End} \\ \hline
\textbf{S1} & Slocum - CTD, Microstructure, ADCP & 2 & 8 & 2025-01-13 & 2025-02-05 \\ \hline
\textbf{S1} & Slocum - CTD, ADCP & 2 & 10 & 2025-01-11 & 2025-02-10 \\ \hline
\textbf{S1} & Slocum - CTD, Optode, PAR, Ecopuck & 2 & 12 & 2025-01-16 & 2025-03-17 \\ \hline
\textbf{S2} & Slocum - CTD, Microstructure, ADCP & 3 & 7 & 2025-01-21 & 2025-03-22 \\ \hline
\textbf{S2} & Slocum - CTD, ADCP & 1 & 9 & 2025-01-31 & 2025-03-22 \\ \hline
\textbf{S2} & Slocum - CTD, Optode, PAR, Ecopuck & 1 & 5 & 2025-01-26 & 2025-03-17 \\ \hline
\textbf{S3} & Slocum - CTD, Microstructure, ADCP & 3 & 10 & 2025-01-16 & 2025-03-22 \\ \hline
\textbf{S3} & Slocum - CTD, ADCP & 3 & 8 & 2025-01-11 & 2025-03-12 \\ \hline
\textbf{S3} & Slocum - CTD, Optode, PAR, Ecopuck & 2 & 7 & 2025-01-11 & 2025-03-17 \\ \hline
\textbf{S1} & ALR1500 - UVP6, EK80, EcoBBRTD & 2 & 20 & 2025-01-11 & 2025-02-20 \\ \hline
\textbf{S2} & ALR1500 - UVP6, EK80, EcoBBRTD & 2 & 15 & 2025-02-10 & 2025-03-02 \\ \hline
\end{tabular}%
}
\caption{Requested science surveys around Rothera}
\label{tab:tasks_case_study1}
\end{table}

\begin{figure}[H]
    \centering
    \includegraphics[width=1\linewidth]{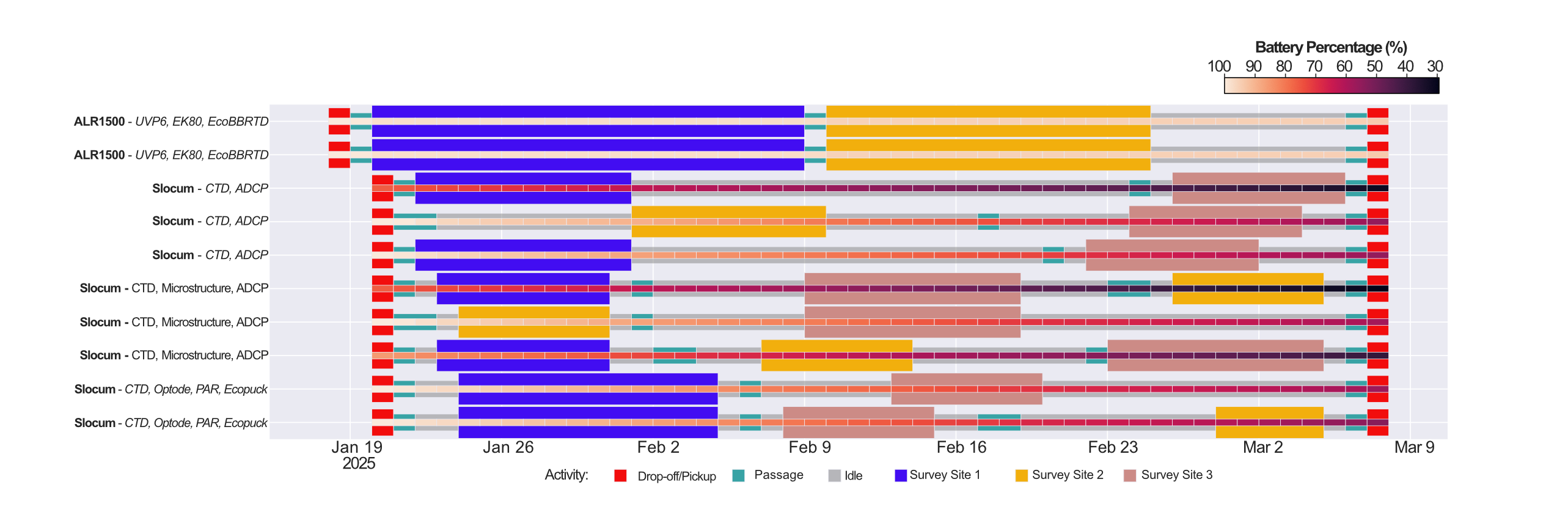}
    \caption{Timeline of MAV activities for the deployment around Rothera}
    \label{fig:timeline_case_study1}
\end{figure}

\begin{sloppypar}
Providing the solver with the science surveys and the fleet information, it determines a fully feasible plan in a matter of seconds, in which every one of the 11 task requirements is satisfied, while also identifying that one MAV of type "\textit{Slocum - CTD+optode+PAR+Ecopuck}" is surplus to requirements and can be safely left ashore. This is a direct, practical benefit of the model: rather than deploying the entire available fleet by default, the optimisation actively searches for a leaner combination of assets, freeing a vehicle for maintenance, a different mission, or simply reducing deployment risk and cost without affecting scientific output.
\end{sloppypar}
Looking beyond task completion, the timeline in Figure \ref{fig:timeline_case_study1} reveals that nearly every MAV experiences periods of idle time between activities. This behaviour is a natural consequence of the interplay between the time windows attached to each task requirement and the comparatively short transit distances around Rothera; a vehicle can often reach its next assignment well before its window opens. Rather than a shortcoming, these idle intervals represent an opportunity: they could be filled with opportunistic science sampling along the way, extracting additional scientific value from a deployment that is already committed and paid for.

\subsection{South Georgia and the Falklands}

\begin{figure}[H]
    \centering
    \includegraphics[width=0.9\linewidth]{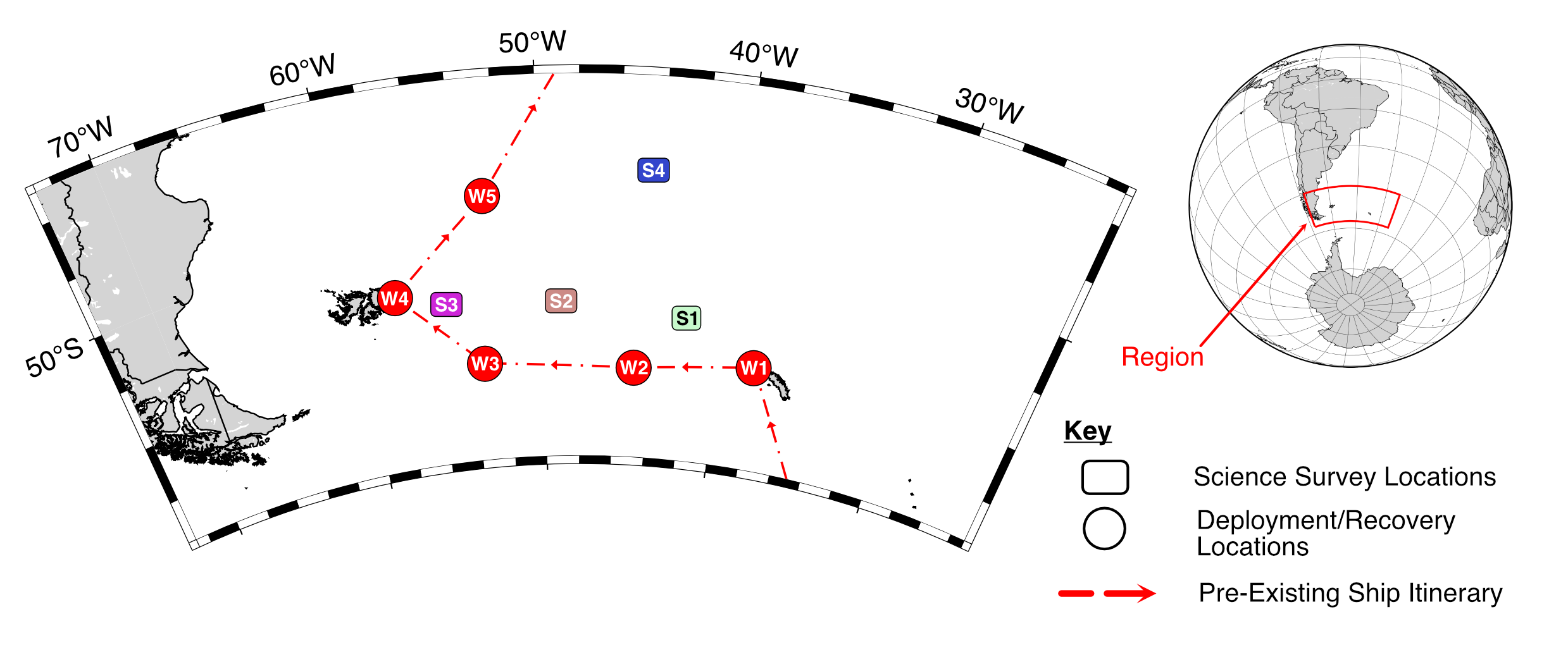}
    \caption{Hybrid ship and shoreside deployment around South Georgia and the Falklands. The red dashed line with arrows represents the pre-defined route of the vessel returning from the Weddell Sea before heading north, travelling through waypoints shown by circles. The research stations are shown by the rectangles.}
    \label{fig:case_study2}
\end{figure}

The second case study demonstrates a hybrid deployment approach, where MAVs can be deployed and recovered from either shoreside locations or the research vessel itself, capturing a more demanding and heterogeneous operational regime than the shoreside-only deployment at Rothera. Distances between locations expand considerably in this scenario, spanning some of the most remote and weather-exposed waters of the Scotia Sea; at this scale, the MAVs' limited battery capacity and the tight time windows attached to each task requirement make direct, unassisted transits infeasible for many routes, so vehicles must instead rendezvous with the ship, ride it as a transit corridor, and resume their science surveys with a freshly swapped battery.

To demonstrate this scenario, we consider a region around South Georgia and the Falklands, where the RRS Sir David Attenborough (SDA) has a pre-defined itinerary returning from the Weddell Sea before heading north (Figure \ref{fig:case_study2}). The starting locations of the 12 available MAVs are spread across: the ship itself (1 Slocum glider), South Georgia (5 Slocum gliders and 2 ALR1500 vehicles) and the Falklands (4 Slocum gliders). The problem consists of 12 task requirements across the 4 research stations and requires a mix of MAV types, with a total of 22 autonomous assets deployed (Table \ref{tab:fleet_case_study2} and Table \ref{tab:tasks_case_study2}); all autonomous assets must rendezvous with the ship before it heads north.

\begin{table}[h]
\resizebox{\textwidth}{!}{%
\begin{tabular}{|cc|c|c|c|}
\hline
\multicolumn{2}{|c|}{\textbf{AUV Types}} & \multirow{2}{*}{\textbf{Base Location}} & \multirow{2}{*}{\textbf{Pickup Location}} & \multirow{2}{*}{\textbf{\begin{tabular}[c]{@{}c@{}}Initial Battery\\ Level (\%)\end{tabular}}} \\ \cline{1-2}
\multicolumn{1}{|c|}{\textbf{Autonomous Platform}} & \textbf{Sensors} &  &  &  \\ \hline
\multicolumn{1}{|c|}{Slocum} & CTD+microstructure+ADCP & Falklands & RDV & 51 \\ \hline
\multicolumn{1}{|c|}{Slocum} & CTD+microstructure+ADCP & Falklands & Checkpoint 1 & 43 \\ \hline
\multicolumn{1}{|c|}{Slocum} & CTD+microstructure+ADCP & South Georgia & Checkpoint 1 & 100 \\ \hline
\multicolumn{1}{|c|}{Slocum} & CTD+microstructure+ADCP & Checkpoint 1 & Checkpoint 1 & 68 \\ \hline
\multicolumn{1}{|c|}{Slocum} & CTD+ADCP & Falklands & RDV & 60 \\ \hline
\multicolumn{1}{|c|}{Slocum} & CTD+ADCP & South Georgia & RDV & 100 \\ \hline
\multicolumn{1}{|c|}{Slocum} & CTD+ADCP & South Georgia & Checkpoint 1 & 72 \\ \hline
\multicolumn{1}{|c|}{Slocum} & CTD+optode+PAR+Ecopuck & Falklands & RDV & 100 \\ \hline
\multicolumn{1}{|c|}{Slocum} & CTD+optode+PAR+Ecopuck & South Georgia & RDV & 67 \\ \hline
\multicolumn{1}{|c|}{Slocum} & CTD+optode+PAR+Ecopuck & Checkpoint 1 & Checkpoint 1 & 54 \\ \hline
\multicolumn{1}{|c|}{ALR1500} & UVP6+EK80+EcoBBRTD & South Georgia & RDV & 100 \\ \hline
\multicolumn{1}{|c|}{ALR1500} & UVP6+EK80+EcoBBRTD & Checkpoint 1 & Checkpoint 1 & 58 \\ \hline
\end{tabular}%
}
\caption{Fleet details for deployment around South Georgia and the Falklands}
\label{tab:fleet_case_study2}
\end{table}

\begin{table}[h]
\resizebox{\textwidth}{!}{%
\begin{tabular}{|c|c|c|c|c|c|}
\hline
\textbf{\begin{tabular}[c]{@{}c@{}}Research \\ Stations\end{tabular}} & \textbf{AUV Types} & \textbf{\begin{tabular}[c]{@{}c@{}}Number \\ Required\end{tabular}} & \textbf{\begin{tabular}[c]{@{}c@{}}Task Duration \\ (days)\end{tabular}} & \textbf{Earliest Start} & \textbf{Latest End} \\ \hline
\textbf{S1} & Slocum - CTD+microstructure+ADCP & 2 & 7 & 2025-01-13 & 2025-08-09 \\ \hline
\textbf{S1} & Slocum - CTD+ADCP & 1 & 10 & 2025-01-16 & 2025-01-27 \\ \hline
\textbf{S1} & Slocum - CTD+optode+PAR+Ecopuck & 1 & 8 & 2025-01-11 & 2025-01-31 \\ \hline
\textbf{S2} & Slocum - CTD+microstructure+ADCP & 2 & 7 & 2025-01-21 & 2025-03-07 \\ \hline
\textbf{S2} & ALR1500 - UVP6+EK80+EcoBBRTD & 1 & 15 & 2025-01-11 & 2025-01-31 \\ \hline
\textbf{S2} & Slocum - CTD+optode+PAR+Ecopuck & 1 & 10 & 2025-01-16 & 2025-01-31 \\ \hline
\textbf{S3} & Slocum - CTD+microstructure+ADCP & 3 & 15 & 2025-05-11 & 2025-06-20 \\ \hline
\textbf{S3} & Slocum - CTD+optode+PAR+Ecopuck & 1 & 10 & 2025-02-10 & 2025-03-02 \\ \hline
\textbf{S3} & ALR1500 - UVP6+EK80+EcoBBRTD & 2 & 20 & 2025-02-10 & 2025-03-02 \\ \hline
\textbf{S4} & Slocum - CTD+microstructure+ADCP & 3 & 7 & 2025-01-21 & 2025-03-07 \\ \hline
\textbf{S4} & Slocum - CTD+ADCP & 2 & 13 & 2025-01-21 & 2025-03-07 \\ \hline
\textbf{S4} & Slocum - CTD+optode+PAR+Ecopuck & 1 & 10 & 2025-01-21 & 2025-03-07 \\ \hline
\end{tabular}%
}
\caption{Requested science surveys around South Georgia and the Falklands}
\label{tab:tasks_case_study2}
\end{table}

Providing the solver with the science surveys and the fleet information, it determines a solution in around 10 seconds. To capture how ship support shapes the deployment plan, the RRS Sir David Attenborough's (SDA) itinerary is overlaid on the timeline alongside the deployed MAVs, with ship-related activities marked separately to highlight the time windows during which vehicles are docked for either charging or transit.

\begin{figure}[H]
    \centering
    \includegraphics[width=1\linewidth]{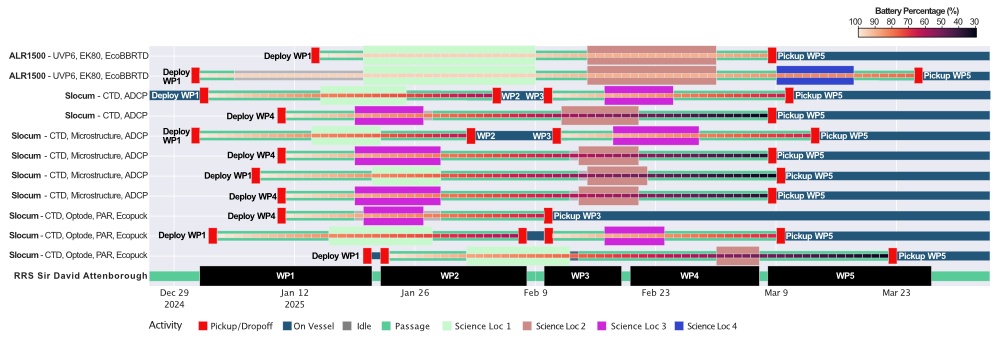}
    \caption{Timeline of MAV activities for the deployment around South Georgia and the Falklands}
    \label{fig:timeline_case_study2}
\end{figure}

The resulting timeline in Figure \ref{fig:timeline_case_study2} shows that the final solution has a complex interaction with the MAVs, where some vehicles leverage the vessel for battery swapping or transit despite being deployed and picked up from shoreside locations. Since the minimum permitted stay aboard the ship is one day, we can distinguish between the two roles the vessel plays: an MAV that boards and disembarks within the same time window is undergoing a battery swap, whereas one that leaves during a later window is instead using the ship purely as a transit shortcut across otherwise unreachable distances.

As with the shoreside deployment at Rothera, the solver also decided that not every available MAV is needed to complete the plan; one Slocum glider is retained aboard the ship rather than deployed, since deploying it would not increase the achieved task coverage.

\begin{figure}[H]
    \centering
    \includegraphics[width=1\linewidth]{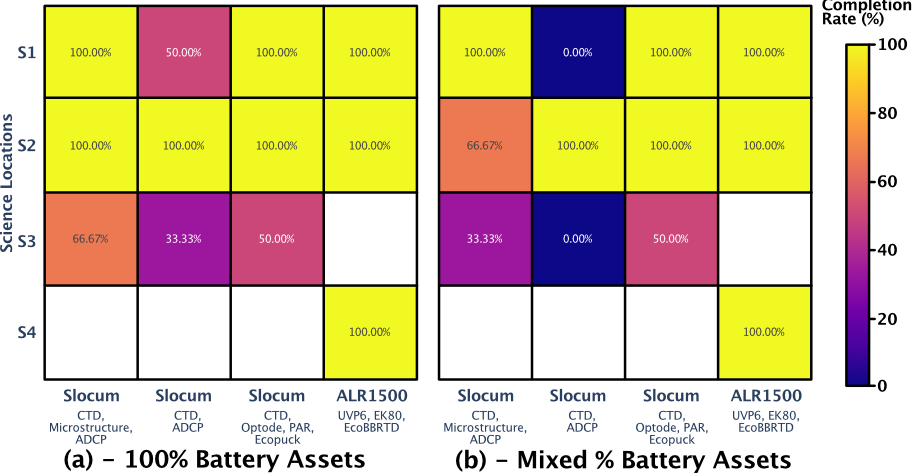}
    \caption{Heatmaps of completion rates of task requirements for the deployment around South Georgia and the Falklands: (a) with the research ship and mixed initial battery levels, and (b) without the research ship and with all MAVs deployed at 100\% of battery}
    \label{fig:heatmap_case_study2}
\end{figure}

\begin{sloppypar}
We generated this example with a fleet consisting of 100\% initial battery levels, but with an infeasible number of science tasks to be completed; even so, the solver can still return a solution that completes as much as possible. The heatmap in Figure \ref{fig:heatmap_case_study2}(a) shows that the task scheduled at research station S1 for MAV type "\textit{Slocum - CTD+microstructure+ADCP}" only reaches half of its required coverage. Because the solution returned is provably optimal given the current inputs, this shortfall is not a modelling artefact but a genuine reflection of the fleet's limits relative to demand; achieving full coverage of this particular task would require either additional assets or a re-weighting of task priorities to favour it over competing requirements.
\end{sloppypar}
To quantify just how much value the ship itself contributes, we re-ran the same scenario with the vessel excluded from the model altogether. To isolate the effect of ship assistance from that of the fleet's initial battery state, we also reset every MAV to a full 100\% charge in this variant, effectively giving the no-ship configuration the best possible starting condition. Despite this advantage, Figure \ref{fig:heatmap_case_study2}(b) shows that overall science coverage still falls short of the original, ship-assisted plan, confirming that the ship is not merely convenient but materially expands the range of achievable science outcomes in this operating environment.

\subsection{United Kingdom Large Scale Deployment}
\label{sssec:uk_large_scale}

\begin{figure}[H]
    \centering
    \includegraphics[width=0.9\linewidth]{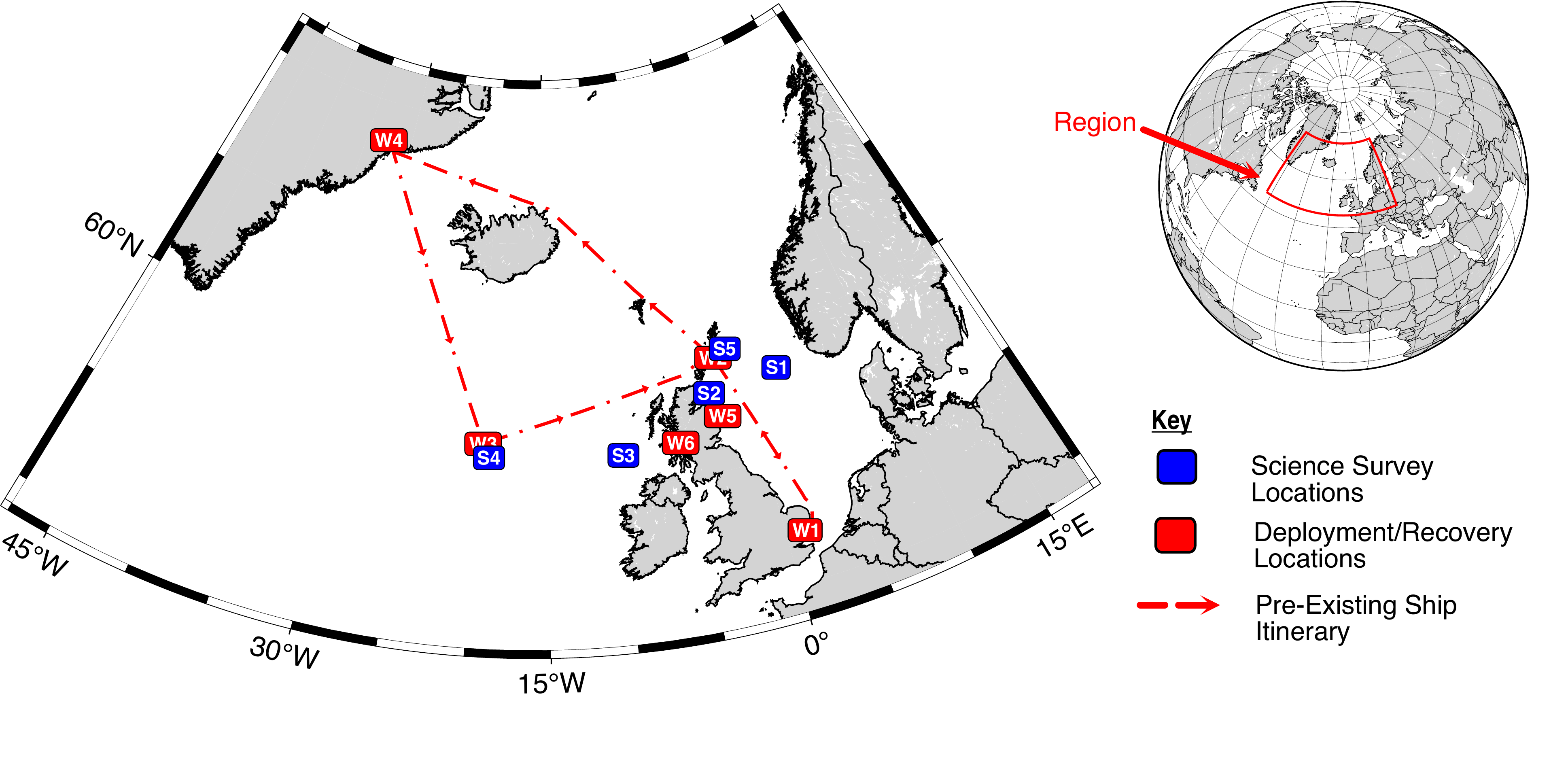}
    \caption{Case study integrating a ship with a pre-defined itinerary. The red dashed line with arrows represents the pre-defined route of the vessel travelling through waypoints shown by circles. The research stations are shown by the rectangles.}
    \label{fig:case_study3}
\end{figure}

The third and last case study demonstrates the scaling of the model to a situation with a large fleet of MAVs (100s) across a variety of start locations and asset ownership, echoing the shared-infrastructure vision set out for the UK's future marine research fleet \citep{parker2021future}. Rather than each organisation planning its own assets in isolation, this scenario asks whether a coordinated, national-scale approach can extract greater scientific value from resources that are traditionally managed and scheduled independently.
We generate a scenario with science surveys distributed around the United Kingdom, involving 100 MAVs drawn from three organisations: the British Antarctic Survey, the National Oceanography Centre (NOC), and the Scottish Association for Marine Science (SAMS). The autonomous fleet consists of 60 Slocum Gliders and 40 ALR1500 vehicles, tasked with completing 17 science surveys (Table \ref{tab:tasks_case_study3}) across 5 research stations (S1-S5). A predefined itinerary for the RRS Sir David Attenborough, departing from and returning to Harwich, threads through a series of science waypoints (Figure \ref{fig:case_study3}), while deployment and pickup of the MAVs can occur either from shoreside locations within Scotland or directly from the vessel, giving the model the opportunity to combine ship-assisted transit with the coordination of a substantially larger and more heterogeneous fleet than in the previous two case studies.

\begin{table}[h]
\resizebox{\textwidth}{!}{%
\begin{tabular}{|c|c|c|c|c|c|}
\hline
\textbf{\begin{tabular}[c]{@{}c@{}}Research \\ Stations\end{tabular}} & \textbf{AUV Types} & \textbf{\begin{tabular}[c]{@{}c@{}}Number \\ Required\end{tabular}} & \textbf{\begin{tabular}[c]{@{}c@{}}Task Duration \\ (days)\end{tabular}} & \textbf{Earliest Start} & \textbf{Latest End} \\ \hline
\textbf{S1} & Slocum - CTD+microstructure+ADCP & 15 & 8 & 2025-01-03 & 2025-02-20 \\ \hline
\textbf{S1} & Slocum - CTD+ADCP & 10 & 10 & 2025-01-16 & 2025-02-20 \\ \hline
\textbf{S1} & Slocum - CTD+optode+PAR+Ecopuck & 8 & 8 & 2025-01-11 & 2025-01-31 \\ \hline
\textbf{S2} & Slocum - CTD+microstructure+ADCP & 7 & 7 & 2025-01-13 & 2025-01-26 \\ \hline
\textbf{S2} & Slocum - CTD+ADCP & 14 & 8 & 2025-01-18 & 2025-02-10 \\ \hline
\textbf{S2} & Slocum - CTD+optode+PAR+Ecopuck & 7 & 10 & 2025-01-16 & 2025-01-31 \\ \hline
\textbf{S3} & Slocum - CTD+microstructure+ADCP & 9 & 15 & 2025-02-05 & 2025-02-25 \\ \hline
\textbf{S3} & Slocum - CTD+optode+PAR+Ecopuck & 10 & 10 & 2025-02-10 & 2025-03-02 \\ \hline
\textbf{S2} & ALR1500 - UVP6+EK80+EcoBBRTD & 30 & 15 & 2025-01-11 & 2025-01-31 \\ \hline
\textbf{S3} & ALR1500 - UVP6+EK80+EcoBBRTD & 22 & 20 & 2025-02-10 & 2025-03-02 \\ \hline
\textbf{S4} & Slocum - CTD+microstructure+ADCP & 16 & 7 & 2025-04-01 & 2025-05-31 \\ \hline
\textbf{S4} & Slocum - CTD+ADCP & 8 & 13 & 2025-03-27 & 2025-05-21 \\ \hline
\textbf{S4} & Slocum - CTD+optode+PAR+Ecopuck & 12 & 10 & 2025-04-06 & 2025-05-26 \\ \hline
\textbf{S5} & Slocum - CTD+ADCP & 14 & 13 & 2025-01-21 & 2025-03-07 \\ \hline
\textbf{S5} & Slocum - CTD+optode+PAR+Ecopuck & 7 & 10 & 2025-01-21 & 2025-03-07 \\ \hline
\textbf{S4} & ALR1500 - UVP6+EK80+EcoBBRTD & 17 & 13 & 2025-02-20 & 2025-04-01 \\ \hline
\textbf{S5} & ALR1500 - UVP6+EK80+EcoBBRTD & 28 & 10 & 2025-01-11 & 2025-03-02 \\ \hline
\end{tabular}%
}
\caption{Requested science surveys around the United Kingdom}
\label{tab:tasks_case_study3}
\end{table}

Providing the solver with the science surveys and the fleet information, it determines a solution in around 10 minutes; the increased solve time compared to the previous two case studies is due to the significant increase in the number of MAVs. The planner is nonetheless still able to return a solution that minimises battery usage and the number of assets while maximising the amount of science delivered, determining that 5 Slocum gliders provide no additional science delivery and should not be used.

\begin{figure}
    \centering
    \includegraphics[width=1\linewidth]{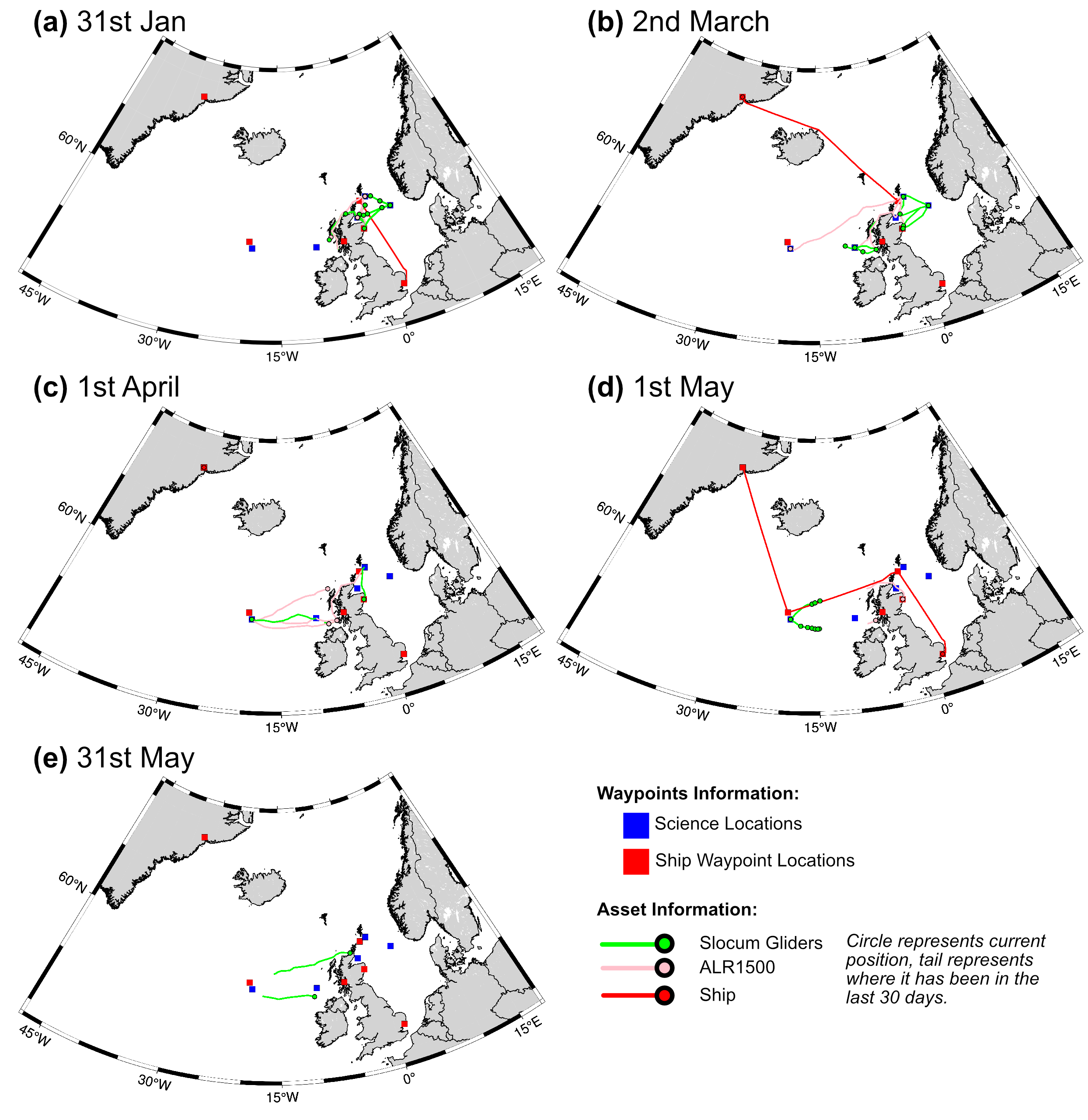}
    \caption{Time evolution of assets in plan. Panels (a) – (e). Each circle represents the asset location on the panel date, with the tails representing its path over the last 30 days, and colours indicating the asset type. The blue boxes represent the science locations and the red boxes the ship waypoint locations.}
    \label{fig:timelocations_case_study3}
\end{figure}

At this scale, with 100 vehicles drawn from three separate organisations, the resulting plan is far too complex to inspect activity-by-activity; instead, we visualise its overall evolution through a series of geographic snapshots of the solver's solution, shown in Figure \ref{fig:timelocations_case_study3}. In each panel, coloured circles mark the position of every asset on a given date, with trailing lines tracing the path each vehicle has followed over the preceding 30 days, allowing the collective behaviour of the fleet to be read at a glance.

The sequence tells a coherent operational story. In the earliest panels, MAVs are concentrated in the North Sea, working through survey tasks close to their point of departure (Panels a and b). As the RRS Sir David Attenborough's itinerary advances, however, the fleet fans out into the more exposed waters of the Atlantic Ocean to complete the remaining surveys, with several vehicles timing their movements to rendezvous with the ship, board for transit, and be redeployed closer to Scotland as the vessel turns back towards Harwich (Panels c and d). The remaining Atlantic-based assets independently make their own way back towards Scotland for recovery (Panels c–e), illustrating how the model orchestrates a mix of independent transits and ship-assisted legs across a large, geographically dispersed fleet without requiring any manual intervention.

\subsection{Worldwide Fleet Coordination}
\label{sssec:worldwide_fleet}

\begin{figure}
    \centering
    \includegraphics[width=0.95\linewidth,height=0.8\textheight,keepaspectratio]{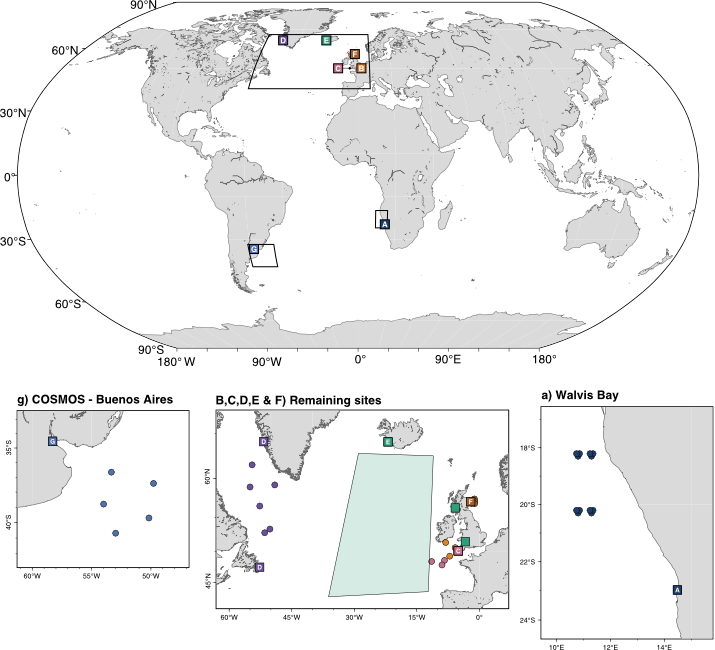}
    \caption{Top: worldwide overview of the seven Case Study 4 sites, lettered A-G and coloured consistently with the zoom panels; marker size is proportional to the number of MAVs deployed at each site, and sites sharing (near-)identical coordinates are spread apart slightly for legibility, connected to their true location by a thin leader line. Bottom: zoomed-in views of each site's deployment/recovery locations (squares) and science survey locations (dots), with sites B-F (all in the North Atlantic/Europe) combined into a single panel.}
    \label{fig:case_study4}
\end{figure}

While the first three case studies each schedule a single fleet at one location, the fourth and final case study demonstrates the applicability of the proposed model across a diverse, worldwide portfolio of MAV missions. Seven independent deployments, each inspired by a real science mission, are planned separately with the proposed model: (i) Walvis Bay (Namibia); (ii) Appledore (Cornwall, UK); (iii) a transit survey between Cornwall and the west coast of France; (iv) REBELS-2 (Nuuk, Greenland); (v) BIO-Carbon ADF (Reykjavik, Iceland); (vi) eSWEETS3 (Aberdeen, Scotland) and (vii) COSMOS (Buenos Aires, Argentina). Figure \ref{fig:case_study4} gives an overview of the geographic spread of the seven deployments (top, marker size proportional to the number of MAVs deployed at each site), together with zoomed-in views of each site's deployment/recovery and science survey locations (bottom): Walvis Bay and COSMOS are shown individually, while the remaining five sites - all clustered around the North Atlantic and north-west Europe - are shown together in a single panel. Each site keeps a consistent colour across both the overview and the zoom panels. This worldwide, multi-fleet scale of coordination is also directly relevant to the national fleet-scaling ambitions set out in Annex A.6 of the UKRI Corporate Plan \citep{UKRI2025CorporatePlan}, which identifies the coordinated growth and utilisation of the UK's national MAV fleet as a strategic priority for future marine research infrastructure.

Table \ref{tab:sites_case_study4} summarises the fleet composition, science task requirements, and duration of each of the seven deployments. Fleet sizes range from a compact 10-vehicle deployment at Appledore and COSMOS up to a 20-vehicle deployment at REBELS-2 and eSWEETS3, while the number of science tasks completed ranges from as few as 6 (Appledore, Cornwall-France) to as many as 92 (BIO-Carbon ADF), reflecting the diversity of mission durations and scientific objectives represented in this case study.

\begin{table}
\resizebox{\textwidth}{!}{%
\begin{tabular}{|c|c|c|c|c|c|c|}
\hline
\textbf{Deployment} & \textbf{Region} & \textbf{Slocum Gliders} & \textbf{ALR1500} & \textbf{Total MAVs} & \textbf{Science Tasks} & \textbf{Total Deployments} \\ \hline
Walvis Bay & Namibia & 8 & 4 & 12 & 36 & 144 \\ \hline
Appledore, Cornwall & Celtic Sea & 6 & 4 & 10 & 6 & 27 \\ \hline
Cornwall to west coast of France & English Channel / Bay of Biscay & 10 & 4 & 14 & 6 & 38 \\ \hline
REBELS-2 (Nuuk, Greenland) & Labrador Sea & 15 & 5 & 20 & 12 & 61 \\ \hline
BIO-Carbon ADF (Reykjavik, Iceland) & North Atlantic & 15 & 4 & 19 & 92 & 98 \\ \hline
eSWEETS3 (Aberdeen, Scotland) & North Sea & 15 & 5 & 20 & 11 & 25 \\ \hline
COSMOS (Buenos Aires, Argentina) & South Atlantic & 5 & 5 & 10 & 10 & 50 \\ \hline
\textbf{Total} & & \textbf{74} & \textbf{31} & \textbf{105} & \textbf{173} & \textbf{443} \\ \hline
\end{tabular}%
}
\caption{Fleet, science task, and deployment summary for the seven Case Study 4 sites}
\label{tab:sites_case_study4}
\end{table}

Across the whole portfolio, the model plans 105 MAVs to complete 173 science tasks over a combined 443 deployments. Fleet utilisation is close to complete: six of the seven deployments use every available asset, with the sole exception of eSWEETS3 (Aberdeen), where 6 Slocum gliders (2 CTD+O2+Fluorescence, 2 CTD+PAR, and 2 ALR1500 UVP6+EK80+EcoBBRTD) are found to be surplus to requirements and safely omitted from the deployed fleet.

\begin{figure}
    \centering
    \includegraphics[width=1.4\linewidth,height=0.8\textheight,keepaspectratio]{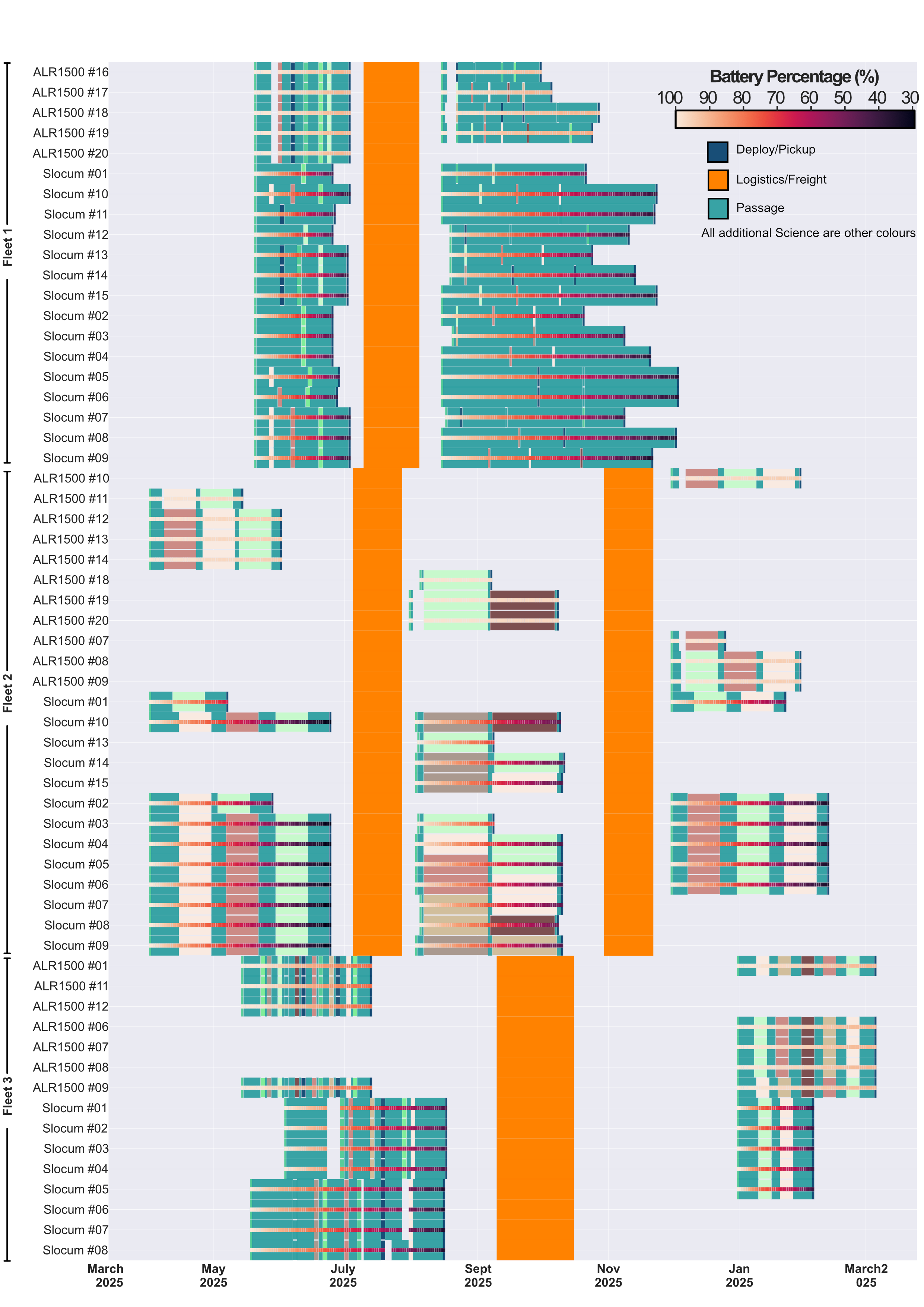}
    \caption{Timeline of activities for the 59 distinct physical vehicles that make up the three fleets serving all seven Case Study 4 deployments (99 asset-deployments in total, with vehicles reused across consecutive projects within their fleet); colour indicates activity type, dark "Logistics" bars mark the transit legs between consecutive deployments of the same fleet, and the thin overlay strip beneath each bar shades the asset's battery level (see colour bar) over time.}
    \label{fig:timeline_case_study4}
\end{figure}

Given the scale of this case study, we first present the timeline of activities across the entire portfolio in a single view. In practice, the seven deployments are not served by seven independent sets of vehicles: the MAVs instead belong to three physical fleets, each of which is shipped or trucked between consecutive deployment sites once one campaign concludes and before the next begins. Fleet 1 covers the North Atlantic, moving by sea from REBELS-2 (Nuuk, Greenland) to BIO-Carbon ADF (Reykjavik, Iceland); Fleet 2 covers the English Channel/North Sea, moving by road between Appledore, the Cornwall-France crossing, and eSWEETS3 (Aberdeen); and Fleet 3 covers the South Atlantic, moving by sea between Walvis Bay (Namibia) and COSMOS (Buenos Aires). Each physical vehicle is identified by its fleet, type, and hull number alone, independent of the specific sensor payload it happened to carry on a given mission, so that the same hull can be recognised as it is redeployed from one science project to the next. Under this view, the 99 actively-deployed MAVs (of the 105 planned; the 6 surplus eSWEETS3 gliders noted above are never activated and so do not appear) reduce to just 59 distinct physical vehicles - 20 in Fleet 1, 24 in Fleet 2, and 15 in Fleet 3 - each reused across an average of nearly two consecutive deployments. Figure \ref{fig:timeline_case_study4} shows the combined schedule of these 59 vehicles, grouped by fleet, with each bar coloured by activity type and shaded by the asset's battery level at that time; the dark "Logistics" bars mark the transit legs connecting consecutive deployments of the same fleet. The plot highlights both the reuse of vehicles within each fleet across sites and the internal structure of each deployment: recurring passage/review cycles as MAVs shuttle between the base and their assigned science survey locations, interspersed with idle periods (grey) where an asset has completed its assigned tasks ahead of the deployment's end date.

This fleet-based logistics structure also points to a substantial carbon saving relative to a conventional single-ship mission, in which a crewed research vessel would itself carry out the science directly - steaming to each station, deploying instrumentation, and recovering samples - rather than merely ferrying a MAV fleet between sites. As a back-of-the-envelope comparison, a mid-sized ocean research vessel of the class needed to run this science programme might reasonably spend around 60\% of the year at sea (roughly 220 days) conducting the surveys itself, burning on the order of 10-12 tonnes of marine fuel per day; combined with a standard marine fuel emission factor of around 3.2 tCO$_2$ per tonne burned \citep{DESNZ2025GHGFactors}, this equates to an annual footprint of roughly 7,500-8,000 tCO$_2$. Under the fleet arrangement modelled here, by contrast, the vessel is only required to freight Fleet 1 and Fleet 3 between consecutive deployment sites (Nuuk to Reykjavik, and Walvis Bay to Buenos Aires, respectively) - a combined transit of around two weeks at sea - while Fleet 2's within-UK relocations are instead carried by road or rail, whichever is the more carbon-efficient option for the route in question and, per standard UK freight conversion factors \citep{DESNZ2025GHGFactors}, over an order of magnitude less carbon-intensive per tonne-kilometre than operating the vessel itself. On this basis, the vessel's own annual footprint falls to roughly 500 tCO$_2$, an overall reduction of the order of 90-95\%, illustrating the scale of carbon saving available simply by decoupling routine ocean data collection from a single, continuously-operating support vessel. This comparison should be read strictly as a proof-of-concept indication of scale rather than a rigorous like-for-like substitution study: a MAV fleet cannot yet replace a research vessel for survey types that depend on sensors or physical sampling equipment the vehicles do not carry - geological core sampling being one prominent example - so in practice a mixed model combining a support vessel with a MAV fleet, rather than full replacement, remains the more realistic near-term outcome.

\begin{figure}
    \centering
    \includegraphics[width=0.85\linewidth]{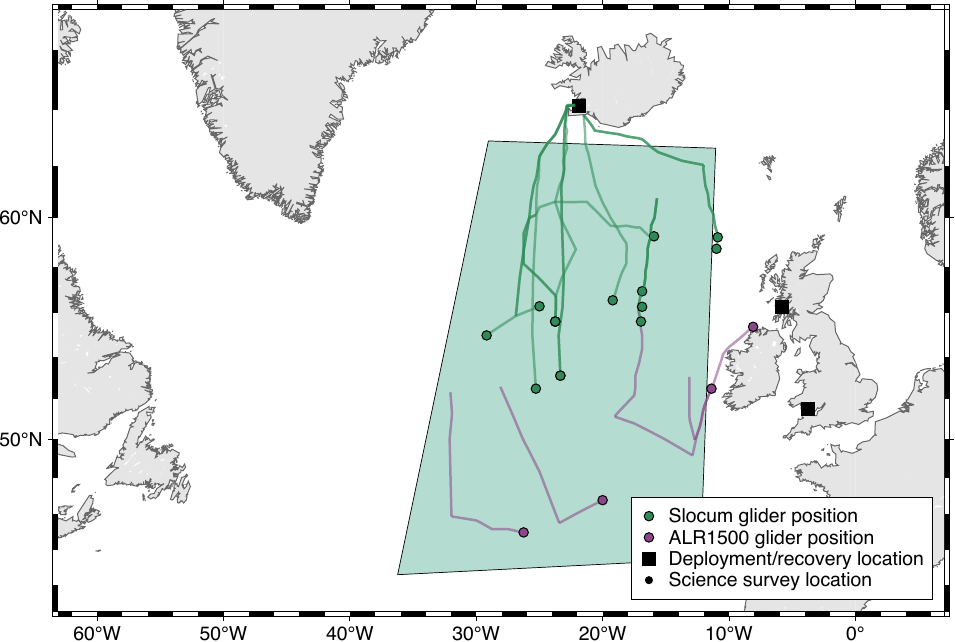}
    \caption{Snapshot of the BIO-Carbon ADF (Reykjavik, Iceland) plan on 15 September 2024, roughly six weeks into the deployment, framed over the same North Atlantic/Europe extent as the combined zoom panel of Figure \ref{fig:case_study4}. Coloured trails show each MAV's PolarRoute-routed track over the preceding two weeks (Slocum gliders in green, ALR1500 vehicles in purple), with a solid marker at its position on the snapshot date; black squares mark the deployment/recovery location and black dots mark science survey locations.}
    \label{fig:stage_case_study4}
\end{figure}

To illustrate the spatial detail available for an individual deployment, Figure \ref{fig:stage_case_study4} shows a single snapshot of the BIO-Carbon ADF (Reykjavik, Iceland) plan - the largest single deployment by number of science tasks - partway through its execution, framed over the same North Atlantic/Europe extent as the combined zoom panel of Figure \ref{fig:case_study4}. Each MAV's most recent two-week track leading up to the snapshot date is drawn as a coloured trail (green for Slocum gliders, purple for ALR1500 vehicles), with a solid marker at its current location; black squares and dots mark the deployment/recovery location and science survey locations respectively. Such snapshots, generated at any date across the full plan, allow the evolving spatial distribution of the fleet to be inspected throughout the deployment; the same visualisation is available for each of the other six deployments in the accompanying case study notebook.

\begin{figure}
    \centering
    \includegraphics[width=0.8\linewidth]{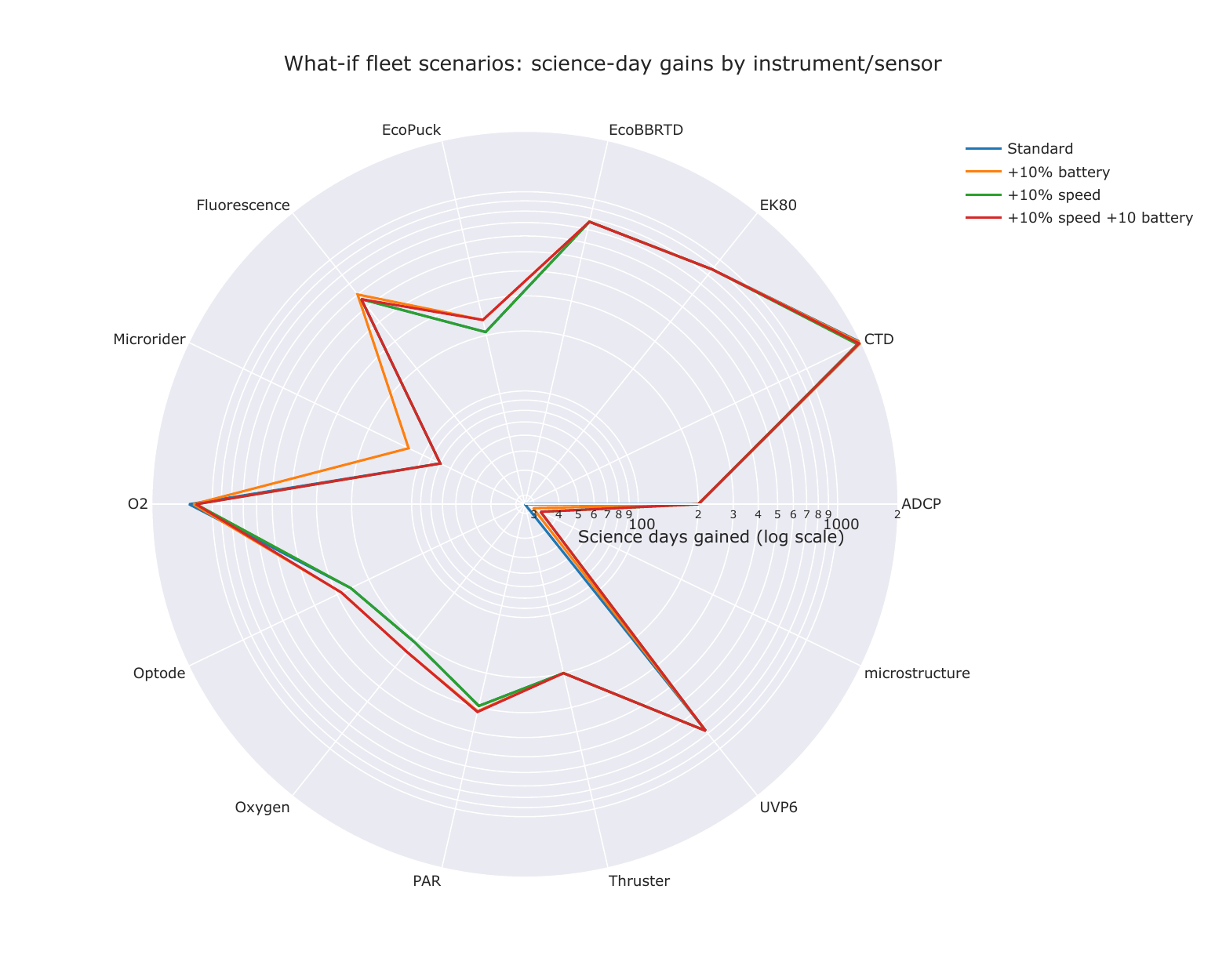}
    \caption{Additional science days gained per instrument/sensor type under three hypothetical fleet upgrades (a 10\% increase in battery capacity, a 10\% increase in transit speed, and both combined), relative to the standard fleet configuration; radial axis on a log scale.}
    \label{fig:radar_case_study4}
\end{figure}

Finally, beyond scheduling a fixed fleet, the same MILP formulation can be used to explore \emph{what-if} fleet design questions - for example, the marginal scientific value of improving a vehicle capability rather than simply adding more vehicles. Figure \ref{fig:radar_case_study4} compares, for one representative deployment, the additional science days gained per instrument/sensor type under three hypothetical fleet upgrades - a 10\% increase in battery capacity, a 10\% increase in transit speed, and both combined - against the standard fleet configuration. Gains are concentrated in a handful of high-value sensors (notably CTD, UVP6, and O2), while several sensors see negligible or no improvement, indicating that they are already limited by factors other than vehicle battery or speed (e.g. science task time windows or site accessibility); this kind of comparison lets a mission planner identify which fleet upgrades would be most cost-effective before committing to procurement.

\section{Discussion and Conclusions}
\label{sec:conclusion}

We proposed in this paper a MILP model for planning marine science missions using a fleet of MAVs. Our model is expected to serve as a decision-support tool for marine operations teams to generate complex plans, as manual planning becomes increasingly expensive. 
The model maximises the number of science tasks covered while optimising the fleet and battery usage. While ship resource constraints are not considered, we allow each MAV to be assigned to the ship only if doing so improves science coverage.
Three case studies inspired by existing marine science missions are solved in this paper. Each one differs from the others in location, science requirements, ship integration, or fleet size. This diversification allows for catching and highlighting the different constraints the model features.

Future improvements in the model could include additional resource constraints, especially for the ship assignments. Similarly, as in the "United Kingdom" case study (Section \ref{sssec:uk_large_scale}), considering multiple ships can also be beneficial in the context of collaborative science missions, where ships and MAVs can be coordinated using a single tool. Since the main objective is maximising scientific gain, and as discussed in Section \ref{sssec:rothera_results}, opportunistic science can benefit the marine science community. Although stakeholders did not raise this, the model can serve as a simulation tool to illustrate how much scientific data we can collect during idle periods.

The case studies presented in this paper also assume a static planning environment, in which the science requirements, environmental conditions and route costs are fixed at the time the plan is generated. In practice, marine science missions face considerable operational uncertainty, and several extensions of the proposed framework would improve its applicability to such dynamic settings. A first extension concerns rapid replanning driven by newly collected data: as MAVs execute a plan, the science data and telemetry they return can reveal that a survey area is of greater or lesser scientific interest than initially assumed, that a task has failed, or that a vehicle's state of charge differs from its forecast, and coupling the model with a data collection update mechanism that feeds these observations back into the planning model and triggers a rapid re-optimisation of the remaining plan would allow the fleet to adapt its behaviour during the mission rather than only at the planning stage. A second extension relates to how routes between locations are represented: the current model assumes a single, precomputed route between each pair of locations, whereas allowing the routing engine to fan out multiple candidate routes for a given leg, rather than committing to one path in advance, would let the planning model select among alternative trajectories according to their cost, risk or robustness, and would provide a natural mechanism for re-routing when conditions change mid-mission. Related to this, route and travel-time costs are currently treated as static inputs, computed once from a fixed environmental snapshot, whereas marine environmental conditions such as currents, wind and sea state vary continuously and are only known through forecasts that are themselves updated over time; extending the model to consume time-varying, forecast-driven edge costs, and to re-solve as updated forecasts become available, would allow plans to reflect the environmental conditions the fleet is actually expected to encounter, rather than a single static assumption. A further extension concerns the representation of science requirements themselves: science tasks are currently modelled as discrete point locations to be visited, yet many marine science missions instead require the survey of an area or transect, in which coverage, rather than the visit of a single point, is the underlying objective, so generalising the task representation to support areal or transect surveys, and the associated notion of partial or ongoing coverage, would broaden the range of missions the model can represent.

Overall, the model and case studies presented in this paper demonstrate that MILP-based planning can generate complex, multi-vehicle marine science plans that would be costly to produce manually, while remaining flexible enough to capture the diverse constraints and priorities of real marine operations. The extensions outlined above, from richer resource and multi-ship modelling to dynamic, forecast-driven and uncertainty-aware replanning, chart a path towards a decision-support tool capable of supporting marine science missions not only at the planning stage but throughout their execution.

\section*{Declaration of AI-assisted tools in the writing process}
During the writing of this paper, the authors used Grammarly for grammar checking and language editing. While using this tool, the authors reviewed the suggested changes to take full responsibility for the publication's content.

\bibliographystyle{apalike}
\bibliography{literature}

\end{document}